\documentclass[11pt]{article}
\usepackage{amsfonts,amssymb,color}
\usepackage[mathscr]{eucal}
\usepackage{amsmath, amsthm}
\usepackage{mathrsfs}
\usepackage{amsbsy}
\usepackage{dsfont}
\usepackage{bbm}
\usepackage{wasysym}
\usepackage{stmaryrd}
\usepackage{upgreek}
\usepackage{hyperref}
\usepackage{comment}
\input xypic
\xyoption{all}
\begin{document}
\baselineskip = 5mm
%%%%%%%%%%%%%%%%%%%%%%%%% fonts %%%%%%%%%%%%%%%%%%%%%%%
\newcommand \ZZ {{\mathbb Z}} % integers
\newcommand \NN {{\mathbb N}} % natural numbers
\newcommand \FF {{\mathbb F}} % finite fields
\newcommand \QQ {{\mathbb Q}} % rational numbers
\newcommand \RR {{\mathbb R}} % real numbers
\newcommand \CC {{\mathbb C}} % complex numbers
\newcommand \PR {{\mathbb P}} % projective space
\newcommand \AF {{\mathbb A}} % affine space
\newcommand \bcA {{\mathscr A}}
\newcommand \bcB {{\mathscr B}}
\newcommand \bcC {{\mathscr C}}
\newcommand \bcF {{\mathscr F}}
\newcommand \bcG {{\mathscr G}}
\newcommand \bcK {{\mathscr K}}
\newcommand \bcN {{\mathscr N}}
\newcommand \bcO {{\mathscr O}}
\newcommand \bcP {{\mathscr P}}
\newcommand \bcR {{\mathscr R}}
\newcommand \bcS {{\mathscr S}}
\newcommand \bcT {{\mathscr T}}
\newcommand \bcU {{\mathscr U}}
\newcommand \bcX {{\mathscr X}}
\newcommand \bcY {{\mathscr Y}}
\newcommand \bcZ {{\mathscr Z}}
\newcommand \catC {{\sf C}}
\newcommand \catD {{\sf D}}
\newcommand \catF {{\sf F}}
\newcommand \catG {{\sf G}}
\newcommand \catE {{\sf E}}
\newcommand \catS {{\sf S}}
\newcommand \catW {{\sf W}}
\newcommand \catX {{\sf X}}
\newcommand \catY {{\sf Y}}
\newcommand \catZ {{\sf Z}}
\newcommand \goa {{\mathfrak a}}
\newcommand \gob {{\mathfrak b}}
\newcommand \goc {{\mathfrak c}}
\newcommand \gom {{\mathfrak m}}
\newcommand \gop {{\mathfrak p}}
\newcommand \goT {{\mathfrak T}}
\newcommand \goC {{\mathfrak C}}
\newcommand \goD {{\mathfrak D}}
\newcommand \goM {{\mathfrak M}}
\newcommand \goN {{\mathfrak N}}
\newcommand \goP {{\mathfrak P}}
\newcommand \goS {{\mathfrak S}}
\newcommand \goH {{\mathfrak H}}
%%%%%%%%%% special symbols %%%%%%%%%%%%%%%%%%%%%%%%%
\newcommand \uno {{\mathbbm 1}}
\newcommand \Le {{\mathbbm L}}
\newcommand \Ta {{\mathbbm T}}
\newcommand \Spec {{\rm {Spec}}}
\newcommand \bSpec {{\bf {Spec}}}
\newcommand \Proj {{\rm {Proj}}}
\newcommand \bProj {{\bf {Proj}}}
\newcommand \Div {{\rm {Div}}}
\newcommand \Pic {{\rm {Pic}}}
\newcommand \Jac {{{J}}}
\newcommand \Alb {{\rm {Alb}}}
\newcommand \NS {{{NS}}}
\newcommand \Corr {{Corr}}
\newcommand \Chow {{\mathscr C}}
\newcommand \Sym {{\rm {Sym}}}
\newcommand \Alt {{\rm {Alt}}}
\newcommand \Prym {{\rm {Prym}}}
\newcommand \cone {{\rm {cone}}}
\newcommand \eq {{\rm {eq}}}
\newcommand \length {{\rm {length}}}
\newcommand \cha {{\rm {char}}}
\newcommand \ord {{\rm {ord}}}
\newcommand \eff {{\rm {eff}}}
\newcommand \shf {{\rm {a}}}
\newcommand \spd {{\rm {s}}}
\newcommand \glue {{\rm {g}}}
\newcommand \equi {{\rm {equi}}}
\newcommand \tr {{\rm {tr}}}
\newcommand \ab {{\rm {ab}}}
\newcommand \add {{\rm {ad}}}
\newcommand \Fix {{\rm {Fix}}}
\newcommand \pty {{\mathbf P}}
\newcommand \type {{\mathbf T}}
\newcommand \prim {{\rm {prim}}}
\newcommand \trp {{\rm {t}}}
\newcommand \cat {{\rm {cat}}}
\newcommand \deop {{\Delta \! }^{op}\, }
\newcommand \pr {{\rm {pr}}}
\newcommand \ev {{\it {ev}}}
\newcommand \defect {{\rm {def}}}
\newcommand \aff {{\rm {aff}}}
\newcommand \Const {{\rm {Const}}}
\newcommand \interior {{\rm {Int}}}
\newcommand \sep {{\rm {sep}}}
\newcommand \td {{\rm {tdeg}}}
\newcommand \tdf {{\mathbf {t}}}
\newcommand \num {{\rm {num}}}
\newcommand \conv {{\it {cv}}}
\newcommand \alg {{\rm {alg}}}
\newcommand \im {{\rm im}}
\newcommand \rat {{\rm rat}}
\newcommand \stalk {{\rm st}}
\newcommand \SG {{\rm SG}}
\newcommand \term {{*}}
\newcommand \Pre {{\mathscr P}}
\newcommand \Funct {{\rm Funct}}
\newcommand \Sets {{\sf Set}}
\newcommand \op {{\rm op}}
\newcommand \Hom {{\rm Hom}}
\newcommand \uHom {{\underline {\rm Hom}}}
\newcommand \HilbF {{\it Hilb}}
\newcommand \HilbS {{\rm Hilb}}
\newcommand \Sch {{\sf Sch}}
\newcommand \cHilb {{\mathscr H\! }{\it ilb}}
\newcommand \cHom {{\mathscr H\! }{\it om}}
\newcommand \cExt {{\mathscr E\! }{\it xt}}
\newcommand \colim {{{\rm colim}\, }}
\newcommand \End {{\rm {End}}}
\newcommand \coker {{\rm {coker}}}
\newcommand \id {{\rm {id}}}
\newcommand \van {{\rm {van}}}
\newcommand \spc {{\rm {sp}}}
\newcommand \Ob {{\rm Ob}}
\newcommand \Aut {{\rm Aut}}
\newcommand \cor {{\rm {cor}}}
\newcommand \res {{\rm {res}}}
\newcommand \tors {{\rm {tors}}}
\newcommand \coeq {{{\rm coeq}\, }}           % coequalizer
\newcommand \Gal {{\rm {Gal}}}
\newcommand \PGL {{\rm {PGL}}}
\newcommand \Gr {{\rm {Gr}}}
\newcommand \Bl {{\rm {Bl}}}
\newcommand \supp {{\rm Supp}}
\newcommand \Sing {{\rm {Sing}}}
\newcommand \spn {{\rm {span}}}
\newcommand \Nm {{\rm {Nm}}}
\newcommand \PShv {{\sf PShv}}                % presheaves
\newcommand \Shv {{\sf Shv}}                  % sheaves
\newcommand \Stk {{\sf Stk}}                  % stacks
\newcommand \sm {{\rm sm}}
\newcommand \reg {{\rm reg}}
\newcommand \nor {{\rm nor}}
\newcommand \noe {{\rm Noe}}
\newcommand \Sm {{\sf Sm}}
\newcommand \Reg {{\sf Reg}}
\newcommand \Nor {{\sf Nor}}
\newcommand \Seminor {{\sf sNor}}
\newcommand \Noe {{\sf Noe}}
\newcommand \inv {{\rm {inv}}}
\newcommand \hc {{\rm {hc}}}
\newcommand \codim {{\rm {codim}}}
\newcommand \ptr {{\pi _2^{\rm tr}}}
\newcommand \Vect {{\mathscr V\! ect}}
\newcommand \ind {{\rm {ind}}}
\newcommand \Ind {{\sf {Ind}}}
\newcommand \Gm {{{\mathbb G}_{\rm m}}}
\newcommand \trdeg {{\rm {tr.deg}}}
\newcommand \seminorm {{\rm {sn}}}
\newcommand \norm {{\rm {norm}}}
\newcommand \Mon {{\sf Mon }}
\newcommand \Mod {{\sf Mod}}
\newcommand \Ab {{\sf Ab }}
\newcommand \tame {\rm {tame }}
\newcommand \prym {\tiny {\Bowtie }}
\newcommand \znak {{\natural }}
\newcommand \et {\rm {\acute e t}}
\newcommand \Zar {\rm {Zar}}
\newcommand \Nis {\rm {Nis}}
\newcommand \Nen {\rm {N\acute en}}
\newcommand \cdh {\rm {cdh}}
\newcommand \h {\rm {h}}
\newcommand \con {\rm {conn}}
\newcommand \sing {{\rm {sing}}}
\newcommand \Top {{\sf {Top}}}
\newcommand \Ringspace {{\sf {Ringspace}}}
\newcommand \qand {{\quad \hbox{and}\quad }}
\newcommand \qqand {{\quad \hbox{and}\quad }}
\newcommand \heither {{\hbox{either}\quad }}
\newcommand \qor {{\quad \hbox{or}\quad }}
\newcommand \Cycl {{\it Cycl }}
\newcommand \PropCycl {{\it PropCycl }}
\newcommand \cycl {{\it cycl }}
\newcommand \PrimeCycl {{\it PrimeCycl }}
\newcommand \PrimePropCycl {{\it PrimePropCycl }}
%%%%%%%%%%%%%%%%%%%%%%%%%%%%%%%%%%%%%%%%%%%%%%%%%%%%
\mathchardef\mhyphen="2D
%%%%%%%%%%%%% theorems, lemmas, etc %%%%%%%%%%%%%%%%
\newtheorem{theorem}{Theorem}
\newtheorem{lemma}[theorem]{Lemma}
\newtheorem{corollary}[theorem]{Corollary}
\newtheorem{proposition}[theorem]{Proposition}
\newtheorem{remark}[theorem]{Remark}
\newtheorem{definition}[theorem]{Definition}
\newtheorem{conjecture}[theorem]{Conjecture}
\newtheorem{example}[theorem]{Example}
\newtheorem{question}[theorem]{Question}
\newtheorem{warning}[theorem]{Warning}
\newtheorem{assumption}[theorem]{Assumption}
\newtheorem{fact}[theorem]{Fact}
\newtheorem{crucialquestion}[theorem]{Crucial Question}
%%%%%%%%%%%%% arrows %%%%%%%%%%%%%%%%%%%%%%
\newcommand \lra {\longrightarrow}
\newcommand \hra {\hookrightarrow}
%%%%%%%%%%%%%%% colors %%%%%%%%%%%
\def\blue {\color{blue}}
\def\red {\color{red}}
\def\green {\color{green}}
%%%%%%%%%%%%%%%%%%%%%%%%%%%%%%%%%%%%%%%%%
\newenvironment{pf}{\par\noindent{\em Proof}.}{\hfill\framebox(6,6)
\par\medskip}
%%%%%%%%%%%%%%%%%%%%%%%%%%%%%%%%%%%%%%%%%%%%%%%%%%%%%%%%%%%%%%%%%%%%%%%%%%%
\title{\bf {The tangent space to the space of 0-cycles}}

\author{Vladimir Guletski\u \i }

\maketitle

\begin{abstract}
\noindent
Let $S$ be a Noetherian scheme, and let $X$ be a scheme over $S$. Under mild assumptions, one can construct the connected infinite symmetric power $\Sym ^{\infty }(X/S)$, whose group completion $\Sym ^{\infty }(X/S)^+$ is an abelian group object in the category of set valued sheaves on the Nisnevich site over $S$. Viewing this completion as the space of relative $0$-cycles on $X/S$, see \cite{SuslinVoevodsky} and \cite{SV-ChowSheaves}, we construct the sheaf of K\"ahler differentials $\Omega ^1_{\Sym ^{\infty }(X/S)^+}$, and the tangent sheaf $T_{\Sym ^{\infty }(X/S)^+}$. We prove that the category of \'etale neighbourhoods at a point $P$ on the space of $0$-cycles is cofiltered. Applying the stalk functor, we also obtain the stalk of the tangent sheaf at $P$, whose tensor product with the residue field is the needed tangent space to the space of $0$-cycles at $P$.
\end{abstract}

%\subjclass[2010]{14A20, 14C25, 14D23, 14J29}

% 14C15     (Equivariant) Chow groups and rings; motives

% 14A20     Generalizations (algebraic spaces, stacks)

% 14C25     Algebraic cycles

% 14D05  	Structure of families (Picard-Lefschetz, monodromy, etc.)

% 14D06  	Fibrations, degenerations

% 14D07  	Variation of Hodge structures [See also 32G20]

% 14D10  	Arithmetic ground fields (finite, local, global)

% 14D15  	Formal methods; deformations [See also 13D10, 14B07, 32Gxx]

% 14D20  	Algebraic moduli problems, moduli of vector bundles For analytic moduli problems, see 32G13

% 14D23 Stacks and moduli problems

% 14E08 Rationality questions

% 14F30 p-adic cohomology, crystalline cohomology

% 14J28 K3-surfaces and Enriques surfaces

% 14J25 Special surfaces

% 14J29 Surfaces of general type

% 14J70 Hypersurfaces

% 14J30 3-folds

% 14J35 4-folds

% 14M20 Rational and unirational varieties

%\keywords{Sheaves, atlases, ringed sites, K\"ahler differentials, tangent sheaf, tangent space, \'etale neighbourhood, cofiltered categories, stalk functor, toposes, locally Noetherian schemes, Nisnevich topology, symmetric powers, free monoids, group completions, relative algebraic cycles, fat points, pullback of relative cycles, relative $0$-cycles, rational equivalence, free rational curve, Bloch's conjecture, surfaces of general type}

\maketitle

\tableofcontents

\section{Introduction}
\label{intro}

The aim of this paper is to make it precise the intuitive feeling that rational equivalence of $0$-cycles on an algebraic variety $X$ is the same as rational connectedness of the corresponding points on a certain space of $0$-cycles on $X$. 

To be more precise, let $X$ be a smooth projective variety over a field $k$. Fix a point on $X$ and use it to embed the $d$-th symmetric power in to the $(d+1)$-th symmetric power of $X$. Passing to a colimit, we obtain the infinite connective symmetric power $\Sym ^{\infty }(X)$ of the variety $X$ over $k$. Viewing it as a commutative monoid, we construct its group completion $\Sym ^{\infty }(X)^+$ in the category of groups. If now $Z_0$ and $Z_1$ are two $0$-cycles on $X$, they can be also considered as points on the group completed symmetric power $\Sym ^{\infty }(X)^+$. Thanks to Suslin-Voevodsky's representability theorem, which holds in arbitrary characteristic, see Theorem 6.8 on page 82 in \cite{SuslinVoevodsky}, $Z_0$ is rationally equivalent to $Z_1$ on $X$ if and only if one can draw a rational curve through $Z_0$ and $Z_1$ on $\Sym ^{\infty }(X)^+$.

This philosophy tracks back, through the paper by Mumford, \cite{Mumford}, to Francesco Severi and possibly earlier. However, it does not give us much, as the object $\Sym ^{\infty }(X)^+$ is not a variety, and it is not clear what could be a rational curve on it, and, more importantly, an appropriate deformation theory of rational curves on the object $\Sym ^{\infty }(X)^+$, in the style of Koll\'ar's approach in \cite{KollarRatCurvesOnVar}. Though Roitman had managed working with the group $\Sym ^{\infty }(X)^+$ as a geometrical object, see \cite{Roitman1} and \cite{Roitman2}, his approach seems to be a compromise, as the necessary technique to deform weird objects was not developed at that time.

Thus, our aim is to develop a technical foundation of deformation theory of rational curves on $\Sym ^{\infty }(X)^+$, and now we are going to explain and justify the concepts promoted in the paper. 

In a broad sense, a geometrical object is a locally ringed site whose Grothendieck topology is of geometric nature. On the other hand, whereas the monoid $\Sym ^{\infty }(X)$ is an ind-scheme, and hence it still can be managed in terms of schemes, the group completion $\Sym ^{\infty }(X)^+$ clearly requires a spacewalk in the category of sheaves on schemes with an appropriate topology, such as \'etale topology or the Nisnevich one. Therefore, we choose that our initial environment is the category of set valued Nisnevich sheaves on locally Noetherian schemes over a base scheme.

Sheaves on a site are still not geometric enough. To produce geometry on a sheaf $\bcX $ we suggest to use the notion of an {\it atlas}, which roughly means that we have a collection of schemes $X_i$ and morphisms of sheaves $X_i\to \bcX $, such that the induced morphism from the coproduct $\coprod _iX_i$ to $\bcX $ is an effective epimorphism. Sheaves with atlases will be called {\it spaces}. This idea gives us a possibility to speak about orphisms from a scheme to a Nisnevich sheaf $\bcX $, and when it is \'etale, with regard to a given atlas on $\bcX $. A Nisnevich-\'etale site $\bcX _{\Nis \mhyphen \et }$ is then the site whose underlying category is the category of morphisms from schemes to $\bcX $, which are \'etale with regard to the atlas on $\bcX $, and whose topology is the restriction of the Nisnevich topology on schemes.

For the local study, let $P$ be a point on $\bcX $, i.e. a morphism from the spectrum of a field to $\bcX $, and let $\bcN _P$ be the category of \'etale neighbourhoods of the point $P$ on the site $\bcX _{\Nis \mhyphen \et }$. An important thing here is that, if the category $\bcN _P$ is cofiltered, we obtain an honest stalk functor at $P$, which yields the corresponding point of the topos of sheaves on the site $\bcX _{\Nis \mhyphen \et }$. 

If now $\bcO _{\bcX }$ is the sheaf of rings on the site $\bcX _{\Nis \mhyphen \et }$, inherited from the regular functions on schemes, its stalk $\bcO _{\bcX \! ,\, P}$ is a local ring, for each point $P$ on $\bcX $. Then $(\bcX _{\Nis \mhyphen \et },\bcO _{\bcX })$ is a locally ringed site. The standard procedure then gives us the sheaf of K\"ahler differentials $\Omega ^1_{\bcX /S}$ and its dual, the tangent sheaf $T_{\bcX /S}$ to the space $\bcX $. 

Applying the stalk at $P$ functor to the sheaf $T_{\bcX /S}$, we obtain the stalk $T_{\bcX \! ,\, P}$, and tensoring by the residue field $\kappa (P)$ of the local ring $\bcO _{\bcX \! ,\, P}$ we obtain the tangent space
  $$
  T_{\bcX }(P)=T_{\bcX \! ,\, P}\otimes \kappa (P)
  $$
to the space $\bcX $ at $P$, with regard to the atlas on $\bcX $. 

Thus, a geometrical object to us is a sheaf $\bcX $ with an atlas, such that $\bcN _P$ is cofiltered for each point $P$ on $\bcX $, and hence the site $\bcX _{\Nis \mhyphen \et }$ is locally ringed by the ring $\bcO _{\bcX }$.

Now, let $S$ be a base scheme, and let $X$ be a locally Noetherian scheme over $S$, such that the relative symmetric power $\Sym ^d(X/S)$ exists for each $d$ (this is always the case if, say, $X$ is quasi-affine or quasi-projective over $S$). 

Assume, moreover, that the structural morphism from $X$ to $S$ admits a section, and use it to construct the monoid $\Sym ^{\infty }(X/S)$, which is and an ind-scheme over $S$. Let $\Sym ^{\infty }(X/S)^+$ be the group completion in the category of Nisnevich sheaves on locally Noetherian schemes over $S$. If $S$ is either of pure characteristic $0$ or $X$ is flat over $S$, then $\Sym ^{\infty }(X/S)^+$ is isomorphic to the sheaf of relative $0$-cycles in the sense of Rydh, see \cite{RydhThesis}. If, moreover, $S$ is seminormal over $\Spec (\QQ )$, then the restriction of the sheaf $\Sym ^{\infty }(X/S)^+$ on schemes seminormal over $S$ gives us a sheaf isomorphic to the sheaves of relative $0$-cycles constructed by Suslin and Voevodsky in \cite{SV-ChowSheaves}, and by Koll\'ar in \cite{KollarRatCurvesOnVar}. 

The fibred squares of symmetric powers yield a natural atlas, we call it the Chow atlas, on the sheaf $\Sym ^{\infty }(X/S)^+$. In the paper, we prove the following result. 

%The main result in the paper is the following theorem.

\medskip

\begin{itemize}

\item[]{}

{\sc Theorem.} {\it Let $X$ be a locally Noetherian scheme over $S$, such that symmetric powers exist, and assume that the structural morphism from $X$ to $S$ admits a section. Let $P$ be a point on the space of $0$-cycles $\Sym ^{\infty }(X/S)^+$. Then the category $\bcN _P$ of \'etale neigbourhoods at $P$ is cofiltered.}

\end{itemize}

\medskip

As a consequence of this theorem, we obtain the tangent space
  $$
  T_{\Sym ^{\infty }(X/S)^+}(P)
  $$
to the space $\Sym ^{\infty }(X/S)^+$ at $P$, according to the procedure described above. 

Assume for simplicity that $S$ is the spectrum of an algebraically closed field $k$, such as $\CC $, $\bar \QQ $ or $\bar \FF _p$, and let $X$ be a smooth projective surface of general type with trivial transcendental part in the second \'etale $l$-adic cohomology group $H^2_{\et }(X,\QQ _l)$. Bloch's conjecture predicts that any two $0$-cycles of the same degree on $X$ are rationally equivalent. Reformulating, the space of $0$-cycles $\Sym ^{\infty }(X)^+$ is rationally connected. 

One way of proving that a variety is rationally connected is that we first find a rational curve on it, and then prove that this curve is sufficiently free on this variety. As we have now tangent spaces at points on the space of $0$-cycles, one can try to do the same on $\Sym ^{\infty }(X)^+$. 

This all should be considered in the context of the ideas developed by Green and Griffiths in \cite{GreenGriffiths}. The problem with their approach is, however, that their tangent space to $0$-cycles is the stalk of a sheaf on the variety itself, but not on a space of $0$-cycles, see, for example, the definition on page 90, or formula (8.1) on page 105 in loc. cit. Moreover, the space of $0$-cycles, as a geometric object, is missing in \cite{GreenGriffiths}. 

Our standpoint is that the space of $0$-cycles should possess explicit 'ontology', to provide a basis to develop analysis on $\Sym ^{\infty }(X)^+$ in the Kock-Lawvere style. 

\bigskip

{\sc Acknowledgements.} I am grateful to Lucas das Dores who spotted a few omissions in the first version of the manuscript.

\bigskip

\section{K\"ahler differentials and spaces of $0$-cycles}
\label{kaehler}

\subsection{Atlases, differentials and tangent spaces}
\label{atdiffandtan}

Throughout the paper we will systematically choose and fix Grothendieck universes, and then working with categories small with regard to these universes, but not mentioning this in the text explicitly. A discussion of the foundational aspects of category theory can be found, for example, in \cite{Shulman} or \cite{therisingsea}.

Let $\catS $ be a topos, and let $\catC $ be a full subcategory in $\catS $, which is closed under finite fibred products. For the purposes which will be clear later, objects in the smaller category $\catC $ will be denoted by Latin letters $X$, $Y$, $Z$ etc, whereas objects in the topos $\catS $ will be denoted by the calligraphic letters, such as $\bcX $, $\bcY $, $\bcZ $ etc.

Let $\tau $ be a topology on $\catC $, and let $\bcO $ be a sheaf of commutative rings on the site $\catC _{\tau }$, which will be considered as the structural sheaf of the ringed site $\catC _{\tau }$. Then $\bcO $ is an object of the topos $\Shv (\catC _{\tau })$ of set valued sheaves on $\catC _{\tau }$, so that the latter is a ringed topos with the structural sheaf $\bcO $.

Given an object $\bcX $ in $\catS $ consider the category $\catC /\bcX $ whose objects are morphisms $X\to \bcX $ in $\catS $, where $X$ are objects of $\catC $, and morphisms are morphism $f:X\to Y$ in $\catC $ over the object $\bcX $. Let $(\catC /\bcX )_{\tau }$ be the big site whose underlying category is $\catC /\bcX $ and the topology on $\catC /\bcX $ is induced by the topology $\tau $ on $\catC $. For short of notation, we denote this site by $\bcX _{\tau }$.

Let also $\bcO _{\bcX }$ be the restriction of the structural sheaf $\bcO $ on the site $\bcX _{\tau }$. We shall look at $\bcO _{\bcX }$ as the structural sheaf of the site $\bcX _{\tau }$. Naturally, $\bcO _{\bcX }$ is an object of the topos $\Shv (\bcX _{\tau })$.

The following definitions are slightly extended versions of the definitions in stack theory. An {\it atlas} $A$ on $\bcX $ is a collection of morphisms
  $$
  A=\{ X_i\to \bcX \} _{i\in I}\; ,
  $$
indexed by a set $I$, such that all the objects $X_i$ are objects of the category $\catC $, the induced morphism
  $$
  e_A:\coprod _{i\in I}X_i\to \bcX
  $$
is an epimorphism in $\catS $, and if
  $$
  X\to \bcX
  $$
is in $A$ and
  $$
  X'\to X
  $$
is a morphism in $\catC $, the composition
  $$
  X'\to X\to \bcX
  $$
is again in $A$. The epimorphism $e_A$ will be called the {\it atlas epimorphism} of the atlas $A$.

Notice that since the category $\catS $ is a topos, and in a topos every epimorphism is regular, for any atlas $A$ on an object $\bcX $ in $\catS $ the atlas epimorphism $e_A$ is a regular epimorphism. Moreover, since every topos is a regular category, and in a regular category regular epimorphisms are preserved by pullbacks, every pullback of $e_A$ is again an epimorphism.

If $A$ is an atlas on $\bcX $ and $B$ is a sunset in $A$, such that $B$ is an atlas on $\bcX $, then we will say that $B$ is a {\it subatlas} on $\bcX $. If $A_0$ is a collection of morphisms from objects of $\catC $ whose coproduct gives an epimorphism onto $\bcX $, the set $A$ of all possible precompositions of morphisms from $A_0$ with morphisms from $\catC $ is an atlas on $\bcX $. We will say that $A$ is generated by the collection $A_0$, and write
  $$
  A=\langle A_0\rangle \; .
  $$
If $A$ consists of all morphisms from objects of $\catC $ to $\bcX $, then we will say that the atlas $A$ is {\it complete}. In contrast, if $A$ is generated by $A_0$ and the latter collection consists of one morphism only, then we will be saying that $A$ is a {\it monoatlas} on the object $\bcX $.

Let
  $$
  f:\bcX \to \bcY
  $$
be a morphism in $\catS $, and assume that the object $\bcY $ has an atlas $B$ on it. We will be saying that $f$ is {\it representable}, with regard to the atlas $B$, if for any morphism
  $$
  Y\to \bcY
  $$
from $B$ the fibred product
  $$
  \bcX \times _{\bcY }Y
  $$
is an object in $\catC $.

Let $\pty $ be a property of morphisms in $\catC $ which is $\tau $-local on the source and target, with regard to the topology $\tau $ and in the sense of Tags 02KO and 036G in \cite{StacksProject}. We will say that the morphism $f:\bcX \to \bcY $ possesses the property $\pty $, with regard to the atlas $B$ on $\bcY $, if (i) $f$ is representable with regard to $B$, and (ii) for any morphism $Y\to \bcY $ from $B$ the base change
  $$
  \bcX \times _{\bcY }Y\to Y
  $$
possesses $\pty $. The stability of $\pty $ under base change and compositions is then straightforward.

Let $\bcX $ and $\bcY $ be objects in $\catS $ and assume that $\bcX $ is endowed with an atlas $A$ and $\bcY $ with an atlas $B$ on them. In such a case the product $\bcX \times \bcY $ also admits an atlas $A\times B$ which consists of products of morphisms from the atlases on $\bcX $ and $\bcY $. We will say the $A\times B$ is the {\it product atlas} on $\bcX \times \bcY $.

For example, if $\bcX $ admits an atlas $A$, the product $\bcX \times \bcX $ admits the square $A\times A$ of the atlas $A$, which is an atlas on $\bcX \times \bcX $. For short, we will write $A^2$ instead of $A\times A$. The diagonal morphism
  $$
  \Delta :\bcX \to \bcX \times \bcX
  $$
is representable with regard to $A^2$ if and only if for any two morphisms
  $$
  X\to \bcX \qqand Y\to \bcX
  $$
from $A$ the fibred product
  $$
  X\times _{\bcX }Y
  $$
is an object in $\catC $. In other words, $\Delta $ is representable with regard to $A^2$ if and only if any morphism from $A$ is representable with regard to $A$. If $\Delta $ is representable with regard to $A^2$ then, for short, we will say that $\bcX $ is {\it $\Delta $-representable} with regard to $A$.

Let $\bcX $ be an object in $\catS $ with an atlas $A$ on it. Let $(\catC /\bcX )_{\pty }$ be the subcategory in $\catC /\bcX $ generated by morphisms $X\to \bcX $ which are representable and possess the property $\pty $ with regard to the atlas $A$ on $\bcX $. Since the property $\pty $ is $\tau $-local on the source and target, the subcategory $(\catC /\bcX )_{\pty }$ is closed under fibred products, and therefore we can restrict the topology $\tau $ from $\catC /\bcX $ to $(\catC /\bcX )_{\pty }$ to obtain a small site $\bcX _{\tau \mhyphen \pty }$. This site depends on the atlas on $\bcX $.

The site $\bcX _{\tau \mhyphen \pty }$ can be further tuned as follows. Let $\type $ be a type of objects in $\catC $, and let $\catC _{\type }$ be the corresponding full subcategory in $\catC $. Assume that $\type $ is closed under fibred products in $\catC $, i.e. for any two morphisms $X\to Z$ and $Y\to Z$ in $\catC _{\type }$ the fibred product $X\times _ZY$ in $\catC $ is again an object of type $\type $. Let $(\catC _{\type }/\bcX )_{\pty }$ be the full subcategory in the category $(\catC /\bcX )_{\pty }$ generated by morphisms $X\to \bcX $ possessing the property $\pty $ and such that $X$ is of type $\type $. Since $\pty $ is $\tau $-local on source and target and type $\type $ is closed under fibred products in $\catC $, the category $(\catC _{\type }/\bcX )_{\pty }$ is closed under fibred products. Then we restrict the topology $\tau $ from the category $(\catC /\bcX )_{\pty }$ to the category $(\catC _{\type }/\bcX )_{\pty }$ and obtain a smaller site $\bcX _{\tau \mhyphen \pty \mhyphen \type }$.

Let $\bcX $ and $\bcY $ be two objects in $\catS $ with atlases $A$ and $B$ respectively, and let
  $$
  f:\bcX \to \bcY
  $$
be a morphism in $\catS $. For any morphism
  $$
  X\to \bcX
  $$
from $\bcX _{\tau \mhyphen \pty \mhyphen \type }$ consider the category
  $$
  X/(\catC _{\type }/\bcY )_{\pty }
  $$
of morphisms
  $$
  X\to Y\to \bcY
  $$
such that the square
  \begin{equation}
  \label{moh}
  \xymatrix{
  X\ar[rr]^-{} \ar[dd]_-{} & & Y \ar[dd]^-{} \\ \\
  \bcX \ar[rr]^-{f} & & \bcY
  }
  \end{equation}
commutes, and the morphism $Y\to \bcY $ is in $\bcY _{\tau \mhyphen \pty \mhyphen \type }$. If the category $X/(\catC _{\type }/\bcY )_{\pty }$ is nonempty, for any morphism $X\to \bcX $ from $\bcX _{\tau \mhyphen \pty \mhyphen \type }$, the morphism $f$ creates a functor
  $$
  f^{-1}:\Shv (\bcY _{\tau \mhyphen \pty \mhyphen \type })\to
  \Shv (\bcX _{\tau \mhyphen \pty \mhyphen \type })
  $$
which associates, to any sheaf $\bcF $ on $\bcY _{\tau \mhyphen \pty \mhyphen \type }$, the sheaf $f^{-1}\bcF $ on $\bcX _{\tau \mhyphen \pty \mhyphen \type }$, such that, by definition
  $$
  f^{-1}\bcG (X\to \bcX )=\colim \bcF (Y\to \bcY )\; ,
  $$
where the colimit is taken over the category $X/(\catC _{\type }/\bcY )_{\pty }$.

If $\bcF $ is a sheaf of rings\footnote{in the paper all rings are commutative rings, if otherwise is not mentioned explicitly} on $\bcY _{\tau }$, it is {\it not} true in general that $f^{-1}\bcF $ is a sheaf of rings on $\bcX _{\tau }$. The reason for that is that the forgetful functor from rings to sets commutes with only filtered colimits, whereas the category $X/(\catC _{\type }/\bcY )_{\pty }$ might be well not filtered. But whenever the category $X/(\catC _{\type }/\bcY )_{\pty }$ is nonempty and filtered, the set $f^{-1}\bcF (X\to \bcX )$ inherits the structure of a ring, and if, moreover, this category is nonempty and filtered for any morphism $X\to \bcX $ from $\bcX _{\tau \mhyphen \pty \mhyphen \type }$ the sheaf $f^{-1}\bcF $ is a sheaf of rings on the site $\bcX _{\tau \mhyphen \pty \mhyphen \type }$.

Let us apply the pullback functor $f^{-1}$ to the structural sheaf of rings $\bcO _{\bcY }$. For each pair of two morphisms
  $$
  X\stackrel{g}{\lra }Y\to \bcY \; ,
  $$
such that the square (\ref{moh}) commutes and the second morphism possesses $\pty $, we have a homomorphism of rings
  $$
  \bcO _{\bcY }(Y\to \bcY )=
  \bcO (Y)\stackrel{\bcO (g)}{\lra }\bcO (X)=
  \bcO _{\bcX }(X\to \bcX )\; .
  $$
Such homomorphisms induce a morphism
  $$
  f^{-1}\bcO _{\bcY }(X\to \bcX )\to
  \bcO _{\bcX }(X\to \bcX )\; ,
  $$
for all morphisms $X\to \bcX $, and hence a morphism of set valued sheaves
  \begin{equation}
  \label{lisichki}
  f^{-1}\bcO _{\bcY }\to \bcO _{\bcX }
  \end{equation}

If we assume that the category $X/(\catC _{\type }/\bcY )_{\pty }$ is nonempty and filtered for every $X\to \bcX $ from $\bcX _{\tau \mhyphen \pty \mhyphen \type }$, the morphism (\ref{lisichki}) is a morphism of ring valued sheaves on the site $\bcX _{\tau \mhyphen \pty \mhyphen \type }$. In such a case, though $f$ does not in general give us a morphism of ring topoi, still we can define the sheaf of K\"ahler differentials on $\bcX _{\tau \mhyphen \pty \mhyphen \type }$ of the morphism $f$ as
  $$
  \Omega ^1_{\bcX /\bcY }=
  \Omega ^1_{\bcO _{\bcX }/f^{-1}\bcO _{\bcY }}\; ,
  $$
in terms of page 115 in the first part of \cite{Illusie} (see also the earlier book \cite{GrothCotang}).

Any Gothendieck topos is a cartesian closed category. In particular, the topos $\Shv (\bcX _{\tau \mhyphen \pty \mhyphen \type })$ is a cartesian closed category, for each object $\bcX $ in $\catS $. The internal Hom-objects are given by the following formula. For any two set valued sheaves $\bcF $ and $\bcG $ on the site $\bcX _{\tau \mhyphen \pty \mhyphen \type }$,
  $$
  \cHom (\bcF ,\bcG )(X\to \bcX )=
  \Hom _{\Shv (\bcX _{\tau \mhyphen \pty \mhyphen \type })}
  (\bcF \times X,\bcG )\; ,
  $$
where $X$ is considered as a sheaf on $\bcX _{\tau \mhyphen \pty \mhyphen \type }$ via the Yoneda embedding. Notice also that, if
  $$
  \Hom _X(\bcF \times X,\bcG \times X)
  $$
is a subset of morphisms from $\bcF $ to $\bcG $ over $X$, i.e. the set of morphisms in the slice category $\Shv (\bcX _{\tau \mhyphen \pty \mhyphen \type })/X$, then
  $$
  \Hom _X(\bcF \times X,\bcG \times X)=
  \Hom _{\Shv (\bcX _{\tau \mhyphen \pty \mhyphen \type })}(\bcF \times X,\bcG )
  $$
for elementary categorical reasons. Then the internal Hom can be equivalently defined by setting
  $$
  \cHom (\bcF ,\bcG )(X\to \bcX )=
  \Hom _X(\bcF \times X,\bcG \times X)\; .
  $$

Now, if the category $X/(\catC _{\type }/\bcY )_{\pty }$ is nonempty and filtered, for every $X\to \bcX $ in $\bcX _{\tau \mhyphen \pty \mhyphen \type }$, so that we have the sheaf of K\"ahler differentials $\Omega ^1_{\bcX /\bcY }$, then we can also define the tangent sheaf on $\bcX _{\tau \mhyphen \pty \mhyphen \type }$ to be the dual sheaf
  $$
  T_{\bcX /\bcY }=
  \cHom (\Omega ^1_{\bcX /\bcY },\bcO _{\bcX })\; .
  $$

If
  $$
  \bcY =Z\in \Ob (\catC _{\type })\; ,
  $$
the category $X/(\catC _{\type }/\bcY )_{\pty }$ has a terminal object
  $$
  \xymatrix{
  X\ar[rr]^-{} \ar[dd]_-{} & & Z \ar[dd]^-{\id } \\ \\
  \bcX \ar[rr]^-{f} & & Z
  }
  $$
And since every category with a terminal object is nonempty and filtered, the morphism (\ref{lisichki}) is a morphism of ring valued sheaves, and we obtain the sheaf of K\"ahler differentials
  $$
  \Omega ^1_{\bcX /Z}\in
  \Ob (\Shv (\bcX _{\tau \mhyphen \pty \mhyphen \type }))
  $$
and the tangent sheaf
  $$
  T_{\bcX /Z}\in
  \Ob (\Shv (\bcX _{\tau \mhyphen \pty \mhyphen \type }))
  $$

The above constructions of K\"ahler differentials and tangent sheaves apply to all kinds of geometric setups, embracing smooth and complex-analytic manifolds in terms of synthetic differential geometry, algebraic varieties, schemes, algebraic spaces, stacks, etc. All we need is to choose an appropriate category $\catC $, a topology $\tau $ on $\catC $, a sheaf of rings $\bcO $ and then take $\catS $ to be the category $\PShv (\catC )$ of set valued presheaves on $\catC $ or, when the topology $\tau $ is subcanonical, the category $\Shv (\catC _{\tau })$ of sheaves on the site $\catC _{\tau }$. If a set valued sheaf $\bcX $ on $\catC _{\tau }$ is endowed with an atlas $A$ of morphisms from objects of the category $\catC $ to $\bcX $, then we will say that $\bcX $ is a {\it space}, with regard to the atlas $A$. In other words, a space to us is a sheaf with a fixed atlas on it.

\subsection{Tangent spaces in \'etale and Nisnevich sites} 

All schemes are separated by default. If $X$ is a scheme and $P$ is a point of $X$ then $\varkappa (P)$ will be the residue field of the scheme $X$ at $P$.

%  \begin{equation}
%  \label{chainofcats1}
%  \Reg \subset \Nor \subset \Noe \subset \Sch \; .
%  \end{equation}

%  \begin{equation}
%  \label{chainofcats2}
%  \Reg /S\subset \Nor /S\subset \Noe /S\subset \Sch /S\; .
%  \end{equation}

%%%  Don't forget that since every regular local ring is  %%%
%%%  integrally closed, every regular scheme is normal.   %%%

Let $\Sch $ be the category of schemes. If $S$ is a scheme, let $\Sch /S$ be the category of schemes over $S$. We will always assume that the base scheme $S$ is Noetherian. Let $\Noe /S$ be the full subcategory in $\Sch /S$ generated by locally Noetherian schemes over $S$. We will also need the full subcategory $\Nor /S$ in $\Noe /S$ generated by locally Noetherian schemes which are locally of finite type over $S$ whose structural morphism is normal in the sense of Tag 0390 in \cite{StacksProject}, the full subcategory $\Reg /S$ in $\Nor /S$ generated by locally Noetherian schemes locally of finite type over $S$ whose structural morphism is regular, in the sense of Tag 07R7 in \cite{StacksProject} (since every regular local ring is integrally closed, every regular scheme is normal). Finally, let $\Sm /S$ be the full subcategory in $\Reg /S$ generated by locally Noetherian schemes locally of finite type over $S$ whose structural morphism is smooth. Recall that every smooth scheme over a field is regular, this is why $\Sm /S$ is indeed a full subcategory in $\Reg /S$. Since every regular scheme over a perfect field is smooth, if the residue fields of points on the base scheme $S$ are perfect, the categories $\Sm /S$ and $\Reg /S$ coincide. Thus, we obtain the following chain of full embeddings
  \begin{equation}
  \label{chainofcats3}
  \Sm /S\subset \Reg /S\subset \Nor /S\subset \Noe /S
  \subset \Sch /S\; .
  \end{equation}

The category $\Sch $ possesses the following well-known topologies: the Zariski topology $\Zar $, $\h $-topology, the \'etale topology $\et $, the Nisnevich topology $\Nis $ and the completely decomposed $\h $-topology denoted by $\cdh $. Notice that only the topologies $\Zar $, $\Nis $ and $\et $ are subcanonical, the topologies $\cdh $ and $\h $ are not subcanonical. The relation between these topologies is given by the chains of inclusions
  \begin{equation}
  \label{chainoftops1}
  \Zar \subset \Nis \subset \et \subset \h
  \end{equation}
and
  \begin{equation}
  \label{chainoftops2}
  \Nis \subset \cdh \subset \h \; .
  \end{equation}

The categories $\Sch /S$ and $\Noe /S$ are obviously closed under fibred products. Moreover, the categories $\Nor /S$, $\Reg /S$ and $\Sm /S$ are also closed under fibred products by Propositions 6.8.2 and 6.8.3 in \cite{EGAIV(2)}. For simplicity of notation, the restrictions of all five topologies from (\ref{chainoftops1}) and (\ref{chainoftops2}) on the categories from (\ref{chainofcats3}) will be denoted by the same symbols.

For our purposes the most convenient setup is this:
  $$
  \catC =\Noe /S\; ,\; \; \tau =\Nis \, ,\; \;
  \pty =\et
  $$
and
  $$
  \type \in
  \{ \sm \, ,\; \reg \, ,\; \nor \, ,\; \noe \} \; ,
  $$
i.e.
  $$
  \catC _{\type }\in
  \{ \Sm /S\, ,\; \Reg /S\, ,\; \Nor /S\, ,\; \Noe /S\} \; .
  $$
Since the Nisnevich topology is subcanonical, we can choose
  $$
  \catS =\Shv ((\Noe /S)_{\Nis })
  $$
to be the category of set valued sheaves on the Nisnevich site $(\Noe /S)_{\Nis }$. If a Nisnevich sheaf $\bcX $ is endowed with an atlas $A$ on it, then we will say that $\bcX $ is a {\it Nisnevich space}, with regard to the atlas $A$. Accordingly, for any Nisnevich space $\bcX $ we have the site
  $$
  \bcX _{\Nis \mhyphen \et \mhyphen \type }
  $$
of morphisms from locally Noetherian schemes of type $\type $ over $S$ to $\bcX $, \'etale with regard to the atlas on $\bcX $, endowed with the induced Nisnevich topology on it.

If
  $$
  \type =\noe \; ,
  $$
i.e.
  $$
  \catC _{\type }=\Noe /S\; ,
  $$
then, for short of notation, we will write
  $$
  \bcX _{\Nis \mhyphen \et }
  $$
for instead of $\bcX _{\Nis \mhyphen \et \mhyphen \noe }$.

Notice also that $S$ is a terminal object in the category $\Noe /S$, and, since any sheaf in $\Shv ((\Noe /S)_{\et })$ is the colimit of representable sheaves, $S$ is also a terminal object in the category $\Shv ((\Noe /S)_{\et })$.

Let $\bcX $ be a Nisnevich sheaf on $\Noe /S$. A {\it point} $P$ on $\bcX $ is an equivalence class of morphisms
  $$
  \Spec (K)\to \bcX
  $$
from spectra of fields to $\bcX $ in the category $\Shv _{\Nis }(\Noe /S)$. Two morphisms
  $$
  \Spec (K)\to \bcX \qqand \Spec (K')\to \bcX
  $$
are said to be equivalent if there exists a third field $K''$, containing the fields $K$ and $K'$, such that the diagram
  $$
  \diagram
  \Spec (K'') \ar[dd]_-{} \ar[rr]^-{} & &
  \Spec (K') \ar[dd]^-{} \\ \\
  \Spec (K) \ar[rr]^-{} & & \bcX
  \enddiagram
  $$
commutes. If a morphism from $\Spec (K)$ to $\bcX $ represents $P$ then, by abuse of notation, we will write
  $$
  P:\Spec (K)\to \bcX \; .
  $$

The set of points on $\bcX $ will be denoted by $|\bcX |$. Certainly, if $\bcX $ is represented by a locally Noetherian scheme $X$ over $S$, then $|\bcX |$ is the set of points of the scheme $X$. A geometric point on $\bcX $ is a morphism from $\Spec (K)$ to $\bcX $, where $K$ is algebraically closed. Any geometric point on $\bcX $ represents a point on $\bcX $, and any point on $\bcX $ is represented by a geometric point.

Fix an atlas $A$ on the sheaf $\bcX $. If a point $P$ on $\bcX $ has a representative
  $$
  \Spec (K)\to \bcX \; ,
  $$
and the latter factors through a morphism from $A$, then we will say that $P$ factors through $A$.

Let $P$ be a point of $\bcX $ which factors through $A$. Choose a representative
  $$
  \Spec (K)\to \bcX
  $$
of the point $P$ with $K$ being algebraically closed. Define a functor
  $$
  u_P:\bcX _{\Nis \mhyphen \et \mhyphen \type }\to \Sets
  $$
sending an \'etale morphism
  $$
  X\to \bcX \; ,
  $$
where $X$ is of type $\type $ over $S$, to the set
  $$
  u_P(X\to \bcX )=|X_P|
  $$
of points on the fibre
  $$
  X_P=X\times _{\bcX }\Spec (K)
  $$
of the morphism $X\to \bcX $ at $P$. Notice that since the morphism $X\to \bcX $ is \'etale, it is representable with regard to the atlas $A$ on the sheaf $\bcX $. And since $P$ factorizes through $A$, the fibre $X_P$ is a locally Noetherian scheme over $S$.

If $X$ and $X'$ are two schemes of type $\type $ over $S$ and endowed with two \'etale morphisms $X\to \bcX $ and $X'\to \bcX $, and if
  $$
  f:X\to X'
  $$
is a morphism of schemes over $S$ and over $\bcX $, i.e. a morphism in $\bcX _{\Nis \mhyphen \et \mhyphen \type }$, then
  $$
  u_P(f):u_P(X)\to u_P(X')
  $$
is the map of sets
  $$
  |X_P|\to |X'_P|
  $$
induced by the scheme-theoretical morphism
  $$
  X_P\to X'_P\; ,
  $$
which is, in turn, induced by the morphism $X\to X'$.

Let $X$ be a locally Noetherian scheme of type $\type $ over $S$, let
  $$
  \{ X_i\to X\} _{i\in I}
  $$
be a Nisnevich covering in $\Noe /S$, and let
  $$
  X\to \bcX
  $$
be a morphism in $\Shv ((\Noe /S)_{\et })$, \'etale with regard to the atlas $A$ on $\bcX $. Since every morphism $X_i\to X$ is smooth, and therefore of type $\type $ over $X$, the cover $\{ X_i\to X\} $ is also a Nisnevich cover of the site $\bcX _{\Nis \mhyphen \et \mhyphen \type }$. Applying the functor $u_P$ we obtain the morphism
  $$
  \coprod _{i\in I}u_P(X_i)\to u_P(X)\; ,
  $$
which is nothing else but the set-theoretical map
  $$
  \coprod _{i\in I}|(X_i)_P|\to |X_P|\; .
  $$
Since $P$ factors through $A$, the latter map is surjective.

If $X'$ is another locally Noetherian scheme of type $\type $ over $S$ and
  $$
  X'\to X
  $$
is a morphism of schemes over $S$, such that the composition
  $$
  X'\to X\to \bcX
  $$
is \'etale with regard to $A$, then we look at the morphism
  $$
  u_P(X_i\times _XX')\to u_P(X_i)\times _{u_P(X)}u_P(X')\; ,
  $$
that is the map
  $$
  |(X_i\times _XX')_P|\to |(X_i)_P|\times _{|X_P|}|X'_P|\; .
  $$
Now again, since $P$ factors through $A$, the latter map is bijective.

In other words, the functor $u_P$ satisfies the items (1) and (2) of Tag 00Y5 in \cite{StacksProject}. The last item (3) of the same definition is satisfied when, for example, the category of neighbourhoods of the point $P$ is cofiltered. Let us discuss item (3) in some more detail.

An \'etale neighbourhood of $P$, in the sense of the site $\bcX _{\Nis \mhyphen \et \mhyphen \type }$, is a pair
  $$
  N=(X\to \bcX ,T\in u_P=|X_P|)\; ,
  $$
where $X$ is of type $\type $ over $S$, $X\to \bcX $ is a morphism over $S$, \'etale with regard to the atlas $A$ on $\bcX $, and $T$ is a point of the scheme $X_P$, represented by, say, the morphism
  $$
  \Spec (\kappa (T))\to X_P\; .
  $$

Equivalently, an \'etale neighbourhood of $P$ is just a commutative diagram of type
  $$
  \diagram
  \Spec (K) \ar[dd]_-{} \ar[rrdd]^-{} & & \\ \\
  X \ar[rr]^-{} & & \bcX
  \enddiagram
  $$
\iffalse %%%%%%%%%%%%%%%%%%%%%%%%%%%%%%%%%%%%%%%%%%%%
  $$
  \diagram
  \Spec (K) \ar[dd]_-{} \ar[rrdd]^-{} & & \\ \\
  X \ar[rr]^-{} \ar[dd]^-{} & & \bcX \ar[lldd]^-{} \\ \\
  S
  \enddiagram
  $$
\fi %%%%%%%%%%%%%%%%%%%%%%%%%%%%%%%%%%%%%%%%%%%%%%%%%
where the morphism $X\to \bcX $ is \'etale with regard to the atlas $A$ on $\bcX $, the morphism $\Spec (K)\to \bcX $ represents the point $P$, and all morphisms are over the base scheme $S$.

If
  $$
  N'=(X'\to \bcX ,T'\in |X_P|)
  $$
is another neighbourhood of $P$, a morphism
  $$
  N\to N'
  $$
is a morphism
  $$
  X\to X'
  $$
over $\bcX $, and hence over $S$, such that, if
  $$
  X_P\to X_{P'}
  $$
is the morphism induced on fibres, the composition
  $$
  \Spec (\kappa (T))\to X_P\to X_{P'}
  $$
represents the point $T'$.

Equivalently, if
  $$
  \diagram
  \Spec (K') \ar[dd]_-{} \ar[rrdd]^-{P} & & \\ \\
  X' \ar[rr]^-{} & & \bcX
  \enddiagram
  $$
is another neighbourhood of $P$, a morphism of neighbourhoods is a morphism
  $$
  X\to X'
  $$
over $\bcX $, and hence over $S$, such that, there is a common field extension $K''$ of $K$ and $K'$, such that $\Spec (K'')\to \bcX $ represents $P$, and the diagram
  $$
  \diagram
  \Spec (K'') \ar[dd]_-{} \ar[rr]_-{}
  & & X' \ar[dd]_-{} \\ \\
  X \ar[rruu]^-{} \ar[rr]^-{} & & \bcX
  \enddiagram
  $$
\iffalse %%%%%%%%%%%%%%%%%%%%%%%%%%%%%%%%%%%%%%%%%%%%%%%%
   $$
   \xymatrix{
   & & \Spec (K'') \ar[lldd]_-{}
   \ar[rrdd]^-{} & & \\ \\
   X \ar[rrdddd]^-{} \ar[rrrr]^-{} \ar[rrdd]_-{} & & & &
   X' \ar[lldddd]^-{} \ar[lldd]^-{} \\ \\
   & & \bcX \ar[dd]^-{} & & \\ \\
   & & S & &
   }
   $$
\fi %%%%%%%%%%%%%%%%%%%%%%%%%%%%%%%%%%%%%%%%%%%%%%%%%%%%%
commutes.

Notice that the above definition of a neighbourhood of a point $P$ on $\bcX $ depends on the functor $u_P$, sending $X\to \bcX $ to $|X_P|$. If we change the functor $u_P$, the notion of neighbourhood will be different, see Tag 00Y3 in \cite{StacksProject}.

Let $\bcN _P$ be the category of neighbourhoods of the point $P$ on $\bcX $, in the sense of the site $\bcX _{\Nis \mhyphen \et \mhyphen \type }$. If $\bcF $ is a set valued sheaf on $\bcX _{\Nis \mhyphen \et \mhyphen \type }$, it is, in particular, a set valued presheaf on the same category, and, as such, it induces a functor
  $$
  \bcF |_{\bcN _P^{\op }}:\bcN _P^{\op }\to \Sets
  $$
sending $N=(X\to \bcX ,T\in |X_P|)$ to $\bcF (X)$ and a morphism $N\to N'$ to the obvious map
  $$
  \bcF (X')\to \bcF (X)\; .
  $$
The stalk functor
  $$
  \stalk _P:
  \Shv (\bcX _{\Nis \mhyphen \et \mhyphen \type })
  \to \Sets
  $$
sends a sheaf $\bcF $ on $\bcX _{\Nis \mhyphen \et \mhyphen \type }$ to the colimit
  $$
  \colim (\bcF |_{\bcN _P^{\op }})
  $$
of the functor $\bcF |_{\bcN _P^{\op }}$.

Once again, we should not forget here that the stack functor $\stalk _P$ depends on the definition of a neighbourhood, and the latter depends on the choice of the functor $u_P$, see Tag 00Y3 in \cite{StacksProject}.

Now, as finite limits commute with filtered colimits, if the category $\bcN _P$ is cofiltered, the stalk functor $\stalk _P$ is left exact, and item (3) of Tag 00Y5 in \cite{StacksProject} holds true as well, and the stalk functor $\stalk _P$ gives rise to a point of the topos $\Shv (\bcX _{\Nis \mhyphen \et \mhyphen \type })$, see Tag 00YA in \cite{StacksProject}. If this is the case, it gives us the well-behaved stalks
  $$
  \bcO _{\bcX ,\, P}=\stalk _P(\bcO _{\bcX })\; ,
  $$
  $$
  \Omega ^1_{\bcX /S,\, P}=\stalk _P(\Omega ^1_{\bcX /S})
  $$
and
  $$
  T_{\bcX /S,\, P}=\stalk _P(T_{\bcX /S})
  $$
at the point $P$.

The latter stalk is not, however, a tangent space to $\bcX $ at $P$. To achieve an honest tangent space we need to observe that, whenever $\bcN _P$ is cofiltered for each $P$, the site $\bcX _{\Nis \mhyphen \et \mhyphen \type }$ is locally ringed in the sense of the definition appearing in Exercise 13.9 on page 512 in \cite{SGA4-1} (see page 313 in the newly typeset version), as well as in the sense of a sightly different Tag 04EU in \cite{StacksProject}. Indeed, any scheme $U$ is a locally ring site with enough points. Applying Tag 04ET in \cite{StacksProject}, we see that for any Zariski open subset $V$ in $U$ and for any function $\bcO _U(V)$ there exists an open covering $V=\cup V_i$ of the set $V$ such that for each index $i$ either $f|_{V_i}$ is invertible or $(1-f)|_{V_i}$ is invertible. If now $U\to \bcX $ is an \'etale morphism from a scheme $U$ to $\bcX $ over $S$, with regard to the atlas on $\bcX $, since
  $$
  \Gamma (U,\bcO _{\bcX })=\Gamma (U,\bcO _U)\; ,
  $$
we obtain item (1) of Tag 04ES in \cite{StacksProject}, and the condition (1) in Tag 04ET is obvious.

Now, since the site $\bcX _{\Nis \mhyphen \et \mhyphen \type }$ is locally ringed, we consider the maximal ideal
  $$
  \gom _{\bcX \! ,\, P}\subset \bcO _{\bcX \! ,\, P}
  $$
and let
  $$
  \kappa (P)=\bcO _{\bcX \! ,\, P}/\gom _{\bcX \! ,\, P}
  $$
be the residue field of the locally ring site at the point $P$. Then we also have two vector spaces
  $$
  \Omega ^1_{\bcX /S}(P)=
  \Omega ^1_{\bcX /S\!,\, P}
  \otimes _{\bcO _P}\kappa (P)
  $$
and
  $$
  T_{\bcX /S}(P)=
  T_{\bcX /S\!,\, P}
  \otimes _{\bcO _P}\kappa (P)
  $$
over the residue field $\kappa (P)$. The latter is our {\it tangent space} to the space $\bcX $ at the point $P$.

\subsection{Relative $0$-cycles after Suslin-Voevodsky and Rydh}
\label{relcycles}

The purpose of this section is discuss various approaches to (pre-)sheaves of relative $0$-cycles on schemes, as they are building blocks of spaces of $0$-cycles. First we will discuss the latest approach presented in \cite{RydhThesis}. Rydh's construction of a sheaf of relative $0$-cycles is compatible with the earlier approaches due to Suslin-Voevodsky, \cite{SV-ChowSheaves}, and Koll\'ar, \cite{KollarRatCurvesOnVar}, if we restrict all the sheaves on seminormal schemes. More details about relative $0$-cycles can be found in \cite{Anderson}.

\bigskip

\noindent {\it Rydh's approach}

\medskip

\noindent So let again $X$ be AF over $S$, and for any nonnegative integer $d$ let $\Gamma ^d(X/S)$ be the $d$-divided power of $X$ over $S$, as explained in Paper I in \cite{RydhThesis}. The infinite coproduct
  $$
  \coprod _{d=0}^{\infty }\Gamma ^d(X/S)
  $$
is a monoid in $\Shv ((\Noe /S)_{\Nis })$.

The canonical morphism
  $$
  (X/S)^d\to \Gamma ^d(X/S)
  $$
is $\Sigma _d$-equivariant on the source, see Prop. 4.1.5 in loc.cit, so that there exists also a canonical morphism
  $$
  \Sym ^d(X/S)\to \Gamma ^d(X/S)\; .
  $$
If the base scheme $S$ is of pure characteristic $0$, or if $X$ is flat over $S$, the latter morphism is an isomorphism of schemes by Corollary 4.2.5 in Paper I in \cite{RydhThesis}. In other words, the divided power $\Gamma ^d(X/S)$ differs from the symmetric power $\Sym ^d(X/S)$ only if the residue fields $\kappa (s)$ can have positive characteristic for points $s\in S$ and, at the same time, $X$ is not flat over $S$. From the point of view of the applications which we have in mind, this is quite a bizarre situations, so that the difference between divided and symmetric powers can be ignored in practice, and we introduce it merely for completeness of the theory.

Now, let $U$ be a locally Noetherian scheme over $S$. According to Paper IV in \cite{RydhThesis}, a {\it relative $0$-cycle of degree $d$} on $X\times _SU$ over $U$ is the equivalence class of ordered pairs $(Z,\alpha )$, where $Z$ is a closed subscheme in $X\times _SU$, such that the composition
  $$
  Z\to X\times _SU\to U
  $$
is a finite, and
  $$
  \alpha :U\to \Gamma ^d(Z/U)
  $$
is a morphism of schemes over $U$. Notice that since the morphism $Z\to U$ is finite, it is AF, and therefore the scheme $\Gamma ^d(Z/U)$ does exist. Two such pairs $(Z_1,\alpha _1)$ and $(Z_2,\alpha _2)$ are said to be equivalent if there is a scheme $Z$ and two closed embeddings $Z\to Z_1$ and $Z\to Z_2$, and a morphism of schemes $\alpha :U\to \Gamma ^d(Z/U)$ over $U$, such that the obvious composition
  $$
  U\stackrel{\alpha }{\lra }\Gamma ^d(Z/U)
  \to \Gamma ^d(Z_i/U)
  $$
is $\alpha _i$ for $i=1,2$, see page 9 in Paper IV in \cite{RydhThesis}. If a relative cycle is represented by a pair $(Z,\alpha )$, we will denote it by $[Z,\alpha ]$.

An important property of divided powers is that if
  $$
  g:U'\to U
  $$
is a morphism of locally Noetherian schemes over $S$, the natural map
  \begin{equation}
  \label{dozhdvlesu}
  \Hom _{U'}(U',\Gamma ^d(X\times _SU/U)\times _UU')\to
  \Hom _{U'}(U',\Gamma ^d(X\times _SU'/U'))
  \end{equation}
is a bijection, see page 12 in paper I in \cite{RydhThesis}. This allows us to define pullbacks of relative $0$-cycles. Indeed, let $[Z,\alpha ]$ be a relative cycle on $X\times _SU$ over $U$. Define $Z'$ and a closed embedding of $Z'$ in to $X\times _SU'$ by the Cartesian square
  $$
  \diagram
  Z' \ar[dd]_-{} \ar[rr]^-{} & & X\times _SU'
  \ar[dd]^-{} \\ \\
  Z \ar[rr]^-{} & & X\times _SU
  \enddiagram
  $$
The composition
  $$
  U'\to U\to \Gamma ^d(Z/U)
  $$
induces the unique morphism
  \begin{equation}
  \label{paporotnik}
  U'\to \Gamma ^d(Z/U)\times _UU'
  \end{equation}
over $U'$ whose composition with the projection onto $\Gamma ^d(Z/U)$ is the initial composition. A particular case of the bijection (\ref{dozhdvlesu}) is the bijection
  \begin{equation}
  \label{dozhdvlesu2}
  \Hom _{U'}(U',\Gamma ^d(Z/U)\times _UU')
  \stackrel{\sim }{\to }
  \Hom _{U'}(U',\Gamma ^d(Z'/U'))
  \end{equation}
Applying (\ref{dozhdvlesu2}) to (\ref{paporotnik}) we obtain the uniquely defined morphism
  $$
  \alpha ':U'\to \Gamma ^d(Z'/U')\; .
  $$
Then
  $$
  g^*[Z,\alpha ]=[Z',\alpha ']
  $$
is, by definition, the pullback of the relative $0$-cycle $[Z,\alpha ]$ along the morphism $g$.

It is easy to verify that such defined pullback is functorial, and we obtain the corresponding set valued presheaf
  $$
  \bcY _{0,d}(X/S):(\Noe /S)^{\op }\to \Sets
  $$
sending any locally Noetherian scheme $U$ over $S$ to the set of all relative $0$-cycles of degree $d$ on $X\times _SU$ over $U$. Let also
  $$
  \bcY _0(X/S)=\coprod _{d=0}^{\infty }\bcY _{0,d}(X/S)
  $$
be the total presheaf of relative $0$-cycles of all degrees.

An important thing here is that the presheaf $\bcY _{0,d}(X/S)$ is represented by the scheme $\Gamma ^d(X/S)$, see Paper I and Paper II in \cite{RydhThesis}. And as the Nisnevich topology is subcanonical, it follows that $\bcY _{0,d}(X/S)$ is a sheaf in Nisnevich topology, i.e. an object of the category $\Shv ((\Noe /S)_{\Nis })$, and the same is true with regard to the presheaf $\bcY _0(X/S)$.

Since each sheaf $\bcY _{0,d}(X/S)$ is represented by the divided power $\Gamma ^d(X/S)$, the sheaf $\bcY _0(X/S)$ is represented by the infinite coproduct $\coprod _{d=0}^{\infty }\Gamma ^d(X/S)$, the sheaf $\bcY _0(X/S)$ is a graded monoid in $\Shv ((\Noe /S)_{\Nis })$, and hence we also have its group completion
  $$
  \bcZ _0(X/S)=\bcY _0(X/S)^+\; .
  $$

Moreover, if the structural morphism $X\to S$ admits a section, the graded monoid $\bcY _0(X/S)$ is pointed, and we can also construct the connective monoid
  $$
  \bcY _0^{\infty }(X/S)=\colim _d\, \bcY _{0,d}(X/S)
  $$
and its group completion
  $$
  \bcZ _0^{\infty }(X/S)=\bcY _0^{\infty }(X/S)^+\; .
  $$

\bigskip

\noindent {\it Suslin-Voevodsky's approach}

\medskip

\noindent For any scheme $X$ let $t(X)$ be the topological space of the scheme $X$, and let $c(X)$ be the set of closed subschemes in $X$. Then we have a map
  $$
  t(X)\to c(X)
  $$
sending any point $\zeta \in X$ to its closure $\overline {\{ \zeta \} }$ with the induced reduced structure of a closed subscheme on it. Let
  $$
  \Cycl ^{\eff }(X)=\NN (t(X))
  $$
be the free monoid generated by points on $X$. Elements of $\Cycl ^{\eff }(X)$ are the {\it effective algebraic cycles}, or simply {\it effective cycles} on the scheme $X$. Let also
  $$
  C^{\eff }(X)=\NN (c(X))
  $$
free monoid generated by closed subschemes of $X$. For any closed subscheme
  $$
  Z\to X\; \in \; c(X)
  $$
let $\zeta _1,\ldots ,\zeta _n$ be the generic points of the irreducible components of the scheme $Z$, let
  $$
  m_i=\length (\bcO _{\zeta _i,Z})
  $$
be the multiplicity of the component $Z_i=\overline {\zeta _i}$ in $Z$, and let
  $$
  \cycl _X(Z)=\sum _im_iZ_i
  $$
be the fundamental class of the closed subscheme $Z$ of the scheme $X$. Then we obtain the standard map
  $$
  \cycl _X:c(X)\to \Cycl ^{\eff }(X)\; ,
  $$
  $$
  Z\mapsto \cycl _X(Z)\; .
  $$
The map $\cycl _X$ extends to the homomorphism of monoids
  $$
  \cycl _X:C^{\eff }(X)\to \Cycl ^{\eff }(X)\; ,
  $$
If
  $$
  C(X)=C^{\eff }(X)^+
  $$
and
  $$
  \Cycl (X)=\Cycl ^{\eff }(X)^+
  $$
then we also have the corresponding homomorphism of abelian groups
  $$
  \cycl _X:C(X)\to \Cycl (X)\; .
  $$

Elements of the free abelian group $\Cycl (X)$ will be called {\it algebraic cycles}, or simply {\it cycles} on the scheme $X$. Points
  $$
  \zeta \in t(X)\; ,
  $$
or, equivalently, their closures
  $$
  Z=\overline {\{ \zeta \} }\; ,
  $$
considered as closed subschemes in $X$ with the induced reduced closed subscheme structure, can be also considered as {\it prime cycles} on $X$. If
  $$
  Z=\sum _im_iZ_i\in \Cycl (X)
  $$
is a cycle on $X$, where $Z_i$ are prime cycles, define its support $\supp (Z)$ to be the union
  $$
  \supp (Z)=\cup _iZ_i\in c(X)
  $$
with the induced reduced structure of a closed subscheme of $X$.

Let $S$ be a Noetherian scheme. A point on $S$ can be understood as a morphism
  $$
  P:\Spec (k)\to S
  $$
from the spectrum of a field $k$ to $S$. A {\it fat point} of $S$ over $P$ is then two morphisms of schemes
  $$
  P_0:\Spec (k)\to \Spec (R)\quad \hbox{and}\quad P_1:\Spec (R)\to S\; ,
  $$
where $R$ is a DVR whose residue field is $k$, such that
  $$
  P_1\circ P_0=P\; ,
  $$
the image of $P_0$ is the closed point of $\Spec (R)$, and $P_1$ sends the generic point $\Spec (R_{(0)})$ to the generic point of the scheme $S$.

Let now
  $$
  f:X\to S
  $$
be a scheme of finite type over $S$, and let
  $$
  Z\to X
  $$
be a closed subscheme in $X$. Let $R$ be a discrete valuation ring,
  $$
  D=\Spec (R)\; ,
  $$
and let
  $$
  g:D\to S
  $$
be a morphism of schemes from $D$ to $S$. Let also
  $$
  \eta =\Spec (R_{(0)})
  $$
be the generic point of $D$,
  $$
  X_D=X\times _SD\; ,\quad Z_D=Z\times _SD\qand
  Z_{\eta }=Z\times _S\eta \; .
  $$
Then there exists a unique closed embedding
  $$
  Z'_D\to Z_D\; ,
  $$
such that its pull-back
  $$
  Z'_{\eta }\to Z_{\eta }
  $$
along the morphism $Z_{\eta }\to Z_D$, is an isomorphism, and the composition
  $$
  Z'_D\to Z_D\to D
  $$
is a flat morphism of schemes, see Proposition 2.8.5 in \cite{EGAIV(2)}.

In particular, one can apply this ``platification" process to a fat point $(P_0,P_1)$ over a point $P\in S$ with $g=P_1$. Let $X_P$ be the fibre of the morphism $X_D\to D$ over the point $P_0$,
  $$
  Z_P=Z_D\times _{X_D}X_P\qand
  Z'_P=Z'_D\times _{Z_D}Z_P\; .
  $$
Since the closed subscheme $Z'_D$ of $X_D$ is flat over $D$, we define the pull-back $(P_0,P_1)^*(Z)$ of the closed subscheme $Z$ to the fibre $X_P$ by the formula
  $$
  (P_0,P_1)^*(Z)=\cycl _{X_P}(Z'_P)\; .
  $$
This gives the definition of a pullback along $(P_0,P_1)$ for primes cycles and, by linearity, extends to a homomorphism
  $$
  (P_0,P_1)^*:\Cycl (X)\to \Cycl (X_P)\; .
  $$

The following definition of Suslin and Voevodsky is of crucial importance, see pp 23 - 24 in \cite{SV-ChowSheaves}.

Let
  $$
  Z=\sum m_iZ_i\in \Cycl (X)
  $$
be a cycle on $X$, and let $\zeta _i$ be the generic point of the prime cycle $Z_i$ for each index $i$. Then $Z$ is said to be a {\it relative cycle} on $X$ over $S$ if:

\begin{itemize}

\item{}
for any generic point $\eta $ of the scheme $S$ there exists $i$, such that
  $$
  f(\zeta _i)=\eta \; ,
  $$

\item{}
for any point $P$ on $S$, and for any two fat points $(P_0,P_1)$ and $(P_0',P_1')$ over $P$,
  $$
  (P_0,P_1)^*(Z)=(P_0',P_1')^*(Z)
  $$
in $\Cycl (X_P)$.

\end{itemize}

The sum of relative cycles is a relative cycle again, and the same for taking the opposite cycle in $\Cycl (X)$. The $0$ in $\Cycl (X)$ is relative by convention. Then we see that relative cycles form a subgroup
  $$
  \Cycl (X/S)=
  \{ Z\in \Cycl (X)\; |\;
  \hbox{$Z$ is relative over $S$} \} \; .
  $$
in $\Cycl (X)$. Let also
  $$
  \Cycl ^{\eff }(X/S)=
  \{ Z=\sum m_iZ_i\in \Cycl (X/S)\; |\;
  m_i\geq 0\; \forall i\; \}
  $$
be a monoid of effective relative cycles in $X$ over $S$.

In general the monoid $\Cycl (X/S)$ is {\it not} a free monoid generated by prime relative cycles, and the group $\Cycl (X/S)$ is {\it not} a free abelian group generated by prime relative cycles.

If $\zeta \in t(X)$, the dimension of $\zeta $ in $X$,
  $$
  \dim (\zeta ,X)\; ,
  $$
is, by definition, the dimension of the closure
  $$
  Z=\overline {\{ \zeta \} }
  $$
inside $X$. A relative cycle
  $$
  Z=\sum m_iZ_i\in \Cycl (X/S)
  $$
is said to be of {\it relative dimension} $r$ if the generic point $\zeta _i$ of each prime cycle $Z_i$ has dimension $r$ in its fibre over $S$. In other words, if
  $$
  \eta _i=f(\zeta _i)\; ,
  $$
we look at the fibre $X_{\eta _i}$ of the morphism $f$ at $\eta _i$. The cycle $Z$ is of relative dimension $r$ over $S$ if
  $$
  \dim (\zeta _i,X_{\eta _i})=r
  $$
for each index $i$. If $Z$ is a relative cycle of relative dimension $r$ on $X$, then we write
  $$
  \dim _S(Z)=r\; .
  $$
Following \cite{SV-ChowSheaves}, p 24, we define
  $$
  \Cycl (X/S,r)=\{ Z\in \Cycl (X/S)\; |\; \dim _S(Z)=r\}
  $$
to be the subset of relative algebraic cycles of relative dimension $r$ on $X$, which is obviously a subgroup in $\Cycl (X/S)$. The definition of
  $$
  \Cycl ^{\eff }(X/S,r)=
  \{ Z=\sum m_iZ_i\in \Cycl (X/S,r)\; |\;
  m_i\geq 0\; \forall i\; \}
  $$
is straightforward.

Notice that if $Z$ is a relative cycle of relative dimension $r$, it does not mean that all the components $Z_i$ are of the same dimension $r$. To pick up equidimensional cycles, we need the following definition. For any point $\zeta \in t(X)$ let
  $$
  \dim (X/S)(x)=\dim _{\zeta }(f^{-1}(f(\zeta )))
  $$
be the dimension of the fibre $f^{-1}(f(\zeta ))$ of the morphism $f$ at $\zeta $. The morphism $f$ is said to be {\it equidimensional} of dimension $r$ if every irreducible component of $X$ dominates an irreducible component of $S$ and the function
  $$
  \dim (X/S):t(X)\to \ZZ
  $$
is constant and equals $r$ for every point $\zeta $ on the scheme $X$. A cycle $Z\in Cycl (X/S)$ is equidimensional of dimension $r$ over $S$ if so is the composition
  $$
  \supp (Z)\to X\to S\; .
  $$
Let then
  $$
  \Cycl _{\equi }(X/S,r)=
  \{ Z\in \Cycl (X/S,r)\; |\;
  \hbox{$Z$ is equidim. of dim. $r$}\} \; .
  $$
Accordingly,
  $$
  \Cycl _{\equi }^{\eff }(X/S,r)=
  \{ Z=\sum _im_iZ_i\in \Cycl _{\equi }(X/S,r)\; |\;
  m_i\geq 0\; \forall i\; \} \; .
  $$
  
\subsection{Relative cycles as spaces of $0$-cycles}
\label{relcycles2}

Next, let
  $$
  U\to S
  $$
be a locally Noetherian scheme over $S$ (not necessarily of finite type over $S$). In \cite{SV-ChowSheaves}, for any cycle
  $$
  Z\in \Cycl (X/S,r)
  $$
Suslin and Voevodsky constructed a uniquely defined cycle
  $$
  Z_U\in \Cycl (X\times _SU/U,r)_{\QQ }\; ,
  $$
a pullback of $Z$ along $U\to S$, such that it is compatible with pullbacks long fat points. Here and below, for any abelian group $A$ we denote by $A_{\QQ }$ the tensor product $A\otimes _{\ZZ }\QQ $.

Thus, following Suslin and Voevodsky, we obtain the obvious presheaf
  $$
  \Cycl (X/S,r)_{\QQ }
  $$
on the category $\Noe /S$, such that for any morphism
  $$
  U\to S
  $$
in $\Noe /S$,
  $$
  \Cycl (X/S,r)_{\QQ }(U)=\Cycl (X\times _SU/U,r)_{\QQ }\; ,
  $$
and the restriction morphisms are induced by the Suslin-Voevodsky's pullbacks of relative cycles.

Following \cite{SV-ChowSheaves}, we will say that the pullback $Z_U$ of a cycle $Z\in \Cycl (X/S,r)$ is {\it integral} if it lies in the image of the canonical homomorphism
  $$
  \Cycl (X\times _SU/U,r)\to
  \Cycl (X\times _SU/U,r)_{\QQ }
  $$
for all schemes $U$ in $\Noe /S$, and define the subgroup
  $$
  z(X/S,r)=\{ Z\in \Cycl (X/S,r)\; |\;
  \hbox{$Z_U$ is integral}\} \; .
  $$
Then $z(X/S,r)$ is an abelian subpresheaf in the presheaf $\Cycl (X/S,r)_{\QQ }$ on the category $\Noe /S$.

Let also
  $$
  z^{\eff }(X/S,r)=
  \{ Z=\sum m_iZ_i\in z(X/S,r)\; |\; m_i\geq 0\; \forall i\}
  $$
and
  $$
  z_{\equi }(X/S,r)=
  \{ Z\in z(X/S,r)\; |\; \hbox{$Z$ is equidim. of dim. $r$ over $S$}\} \; .
  $$
Clearly, $z^{\eff }(X/S,r)$ is a subpresheaf of monoids and $z^{\equi }(X/S,r)$ is a presheaf of abelian groups in $z(X/S,r)$.

For any morphism
  $$
  U\to S\; ,
  $$
which is an object of $\Noe /S$, set
  $$
  \PrimeCycl (X\times _SU/U,r)=
  \{ Z\in \Cycl (X\times _SU/U,r)\; |\;
  \hbox{$Z$ is prime}\}
  $$
and
  $$
  \PrimeCycl _{\equi }(X\times _SU/U,r)=
  \{ Z\in \PrimeCycl (X\times _SU/U,r)\; |\;
  \hbox{$Z$ is equidim.}\}
  $$
If $S$ is regular, and if the morphism $U\to S$ is an object of $\Reg /S$, then
  $$
  z^{\eff }(X/S,r)=
  \NN (\PrimeCycl _{\equi }(X\times _SU/U,r))\; ,
  $$
and
  $$
  z_{\equi }(X/S,r)=
  \NN (\PrimeCycl _{\equi }(X\times _SU/U,r))^+\; ,
  $$
see Corollary 3.4.5 in \cite{SV-ChowSheaves}.

It does not mean, however, that $z^{\eff }(X/S,r)$ is a free monoid in the category of set valued presheaves freely generated by a set valued ``presheaf of relative prime cycles of dimension $r$" on the category $\Reg /S$, as the Suslin-Voevodsky pullback of a relative prime cycle is not necessarily a prime cycle, so that the needed set valued presheaf does not exist. But $z_{\equi }(X/S,r)$ is certainly the group completion of $z^{\eff }(X/S,r)$ as a presheaf on $\Reg /S$.

\begin{theorem}
\label{susvoe_sheaf}
Let $S$ be a Noetherian scheme, and let $X$ be a scheme of finite type over $S$. Then the presheaves $z(X/S,r)$ and $z^{\eff }(X/S,r)$ are sheaves in $\cdh $-topology and, as a consequence, in the Nisnevich topology on the category $\Noe /S$.
\end{theorem}

\begin{pf}
See Theorem 4.2.9(1) on page 65 in \cite{SV-ChowSheaves}.
\end{pf}

Relative cycles can be classified by their degrees, provided there exists a projective embedding of $X$ over $S$. Indeed, assume that $X$ is projective over $S$, i.e. there is a closed embedding
  $$
  i:X\to \PR ^n_S
  $$
over $S$. For each cycle
  $$
  Z=\sum m_jZ_j\in \Cycl (X/S)
  $$
one can define its degree
  $$
  \deg (Z,i)=\sum \deg (i(Z_j))
  $$
with regard to the embedding $i$. Let also
  $$
  z^{\eff }_d((X,i)/S,r)=\{ Z\in z_{\equi }(X/S,r)\; |\; \deg (Z,i)=d\} \; .
  $$
The set valued presheaf
  $$
  z^{\eff }_d((X,i)/S,r):\Noe /S\to \Sets
  $$
is given by the formula
  $$
  z^{\eff }_d((X,i)/S,r)(U)=
  \{ Z\in z_{\equi }(X\times _SU/U,r)\; |\;
  \deg (Z,i\times _S\id _U)=d\} \; ,
  $$
for any locally Noetherian scheme $U$ over $S$.

\medskip

Now recall that if $\bcF $ is a set-valued presheaf on $\Noe /S$ then $\bcF $ is said to be $\h $-representable if there is a scheme $Y$ over $S$, such that the $\h $-sheafification $\bcF _{\h }$ of the sheaf $\bcF $ is isomorphic to the $\h $-sheafification $\Hom _S(-X)_{\h }$ of the representable presheaf $\Hom _S(-X)$, see Definition 4.4.1 in \cite{SV-ChowSheaves}.

\begin{theorem}
\label{hrep}
Let $X$ be a projective scheme of finite type over $S$ and fix a projective embedding $i:X\to \PR ^n_S$ over $S$. Then, for any two nonnegative integers $r$ and $d$, the presheaf $z^{\eff }_d((X,i)/S,r)$ is $\h $-representable by a scheme $C_{r,d}(X/S,i)$ projective over $S$, i.e. there is an isomorphism
  $$
  z^{\eff }_d((X,i)/S,r)_{\h }\simeq
  \Hom _S(-,C_{r,d}(X/S,i))_{\h }
  $$
of set valued sheaves in $\h $-topology on $\Noe /S$. Moreover,
  $$
  z^{\eff }(X/S,r)=
  \coprod _{d=0}^{\infty }z^{\eff }_d((X,i)/S,r)\; ,
  $$
and then $z^{\eff }(X/S,r)$ is $\h $-representable by the scheme
  $$
  C_r(X/S)=\coprod _{d=0}^{\infty }C_{r,d}(X/S,i)\; .
  $$
\end{theorem}

\begin{pf}
See Section 4.2 in \cite{SV-ChowSheaves}.
\end{pf}

A disadvantage of Theorem \ref{hrep} is in the presence of $\h $-sheafification. The latter is a retribution for the generality of the representability result. For relative $0$-cycles this obstacle can be avoided as follows.

Recall that we have already defined the category $\Nor /S$, a full subcategory in $\Sch /S$ generated by schemes over $S$ whose structural morphism is normal, i.e. the fibre at every point is a normal scheme, see Tag 0390 in \cite{StacksProject}. Similarly, one can define the notion of a seminormal morphism and introduce a full subcategory $\Seminor /S$ generated by locally Noetherian schemes over $S$ whose structural morphisms are seminormal, so that we have a chain of subcategories
  $$
  \Nor /S\subset \Seminor /S\subset \Noe /S\; .
  $$
For any presheaf $\bcF $ on $\Noe /S$ let $\bcF |_{\Seminor /S}$ be the restriction of $\bcF $ on the subcategory $\Seminor /S$.

To avoid divided powers, suppose that either the base scheme $S$ is of pure characteristic $0$ or $X$ is flat over $S$. Recall that it follows that
  $$
  \Gamma ^d(X/S)=\Sym ^d(X/S)
  $$
by Corollary 4.2.5 in Paper I in \cite{RydhThesis}, and hence one can work with symmetric powers instead of divided ones. By Theorem 3.1.11 on page 30 of the same paper, we have the canonical identifications
  \begin{equation}
  \label{maincanonical1*}
  \bcY _{0,d}(X/S)=\Sym ^d(X/S)\; ,
  \end{equation}
  $$
  \bcY _0(X/S)=\left( \coprod _{d=0}^{\infty }
  \Sym ^d(X/S)\right)\; ,
  $$
  $$
  \bcY _0^{\infty }(X/S)=\Sym ^{\infty }(X/S)\; ,
  $$
  $$
  \bcZ _0(X/S)=
  \left( \coprod _{d=0}^{\infty }\Sym ^d(X/S)\right)^+
  $$
and
  $$
  \bcZ _0^{\infty }(X/S)=
  \Sym ^{\infty }(X/S)^+\; .
  $$
In other words, we do not need $\h $-sheafification to prove representability of sheaves of $0$-cycles in Rydh's terms.

The point here is that, assuming that $S$ is semi-normal over $\Spec (\QQ )$, after restricting of these five sheaves on the category $\Seminor /S$, we also have the corresponding canonical isomorphisms
  \begin{equation}
  \label{maincanonical1'}
  \bcY _{0,d}(X/S)|_{\Seminor /S}\simeq
  z^{\eff }_d((X,i)/S,0)|_{\Seminor /S}\; ,
  \end{equation}
  \begin{equation}
  \label{maincanonical2'}
  \bcY _0(X/S)|_{\Seminor /S}\simeq
  z^{\eff }(X/S,0)|_{\Seminor /S}\; ,
  \end{equation}
  \begin{equation}
  \label{maincanonical3'}
  \bcY _0^{\infty }(X/S)|_{\Seminor /S}\simeq
  z^{\eff }(X/S,0)_{\infty }|_{\Seminor /S}\; ,
  \end{equation}
  \begin{equation}
  \label{maincanonical3.5'}
  \bcZ _0(X/S)|_{\Seminor /S}\simeq
  z(X/S,0)|_{\Seminor /S}
  \end{equation}
and
  \begin{equation}
  \label{maincanonical4'}
  \bcZ _0^{\infty }(X/S)|_{\Seminor /S}\simeq
  z(X/S,0)_{\infty }|_{\Seminor /S}\; .
  \end{equation}
Moreover, the same result holds true when we compare Rydh's sheaves of $0$-cycles with Koll\'ar's sheaves constructed in Chapter I of the book \cite{KollarRatCurvesOnVar}. These important comparison results are proven in Section 10 of Paper IV in \cite{RydhThesis}.

Thus, since now we will always assume that either the base scheme $S$ is of pure characteristic $0$ or $X$ is flat over $S$, to work with symmetric powers, and in all cases when $S$ will be semi-normal over $\QQ $, we will systematically identify the restrictions of Suslin-Voevodsky's and Rydh's sheaves of $0$-cycles on semi-normal schemes via the isomorphisms (\ref{maincanonical1'}), (\ref{maincanonical2'}), (\ref{maincanonical3'}), (\ref{maincanonical3.5'}) and (\ref{maincanonical4'}).

The Nisnevich sheaf $\Sym ^{\infty }(X/S)^+$ will be now used to construct what then will be the most preferable reincarnation of the space of $0$-cycles on $X$ over the base scheme $S$.

\section{The tangent space to the Nisnevich space of $0$-cycles}

\subsection{Chow atlases on the Nisnevich spaces of $0$-cycles}

To consider the sheaf $\Sym ^{\infty }(X/S)^+$ as a geometrical object, we need to endow it with an atlas, in the line of the definitions in Section \ref{kaehler}. The aim of this section is to present a natural atlas, the Chow atlas, on the sheaf of $0$-cycles $\Sym ^{\infty }(X/S)^+$.

First of all, the sheaf of $0$-cycles possesses a natural inductive structure on it. For each non-negative integer $d$ let
  $$
  \iota _d:\Sym ^d(X/S)\to \Sym ^{\infty }(X/S)
  $$
be the canonical morphism in to the colimit. For short of notation, let also
  $$
  \Sym ^{d,d}(X/S)=
  \Sym ^d(X/S)\times _S\Sym ^d(X/S)\; ,
  $$
  $$
  \Sym ^{\infty ,\infty }(X/S)=
  \Sym ^{\infty }(X/S)\times _S\Sym ^{\infty }(X/S)
  $$
and let
  $$
  \iota _{d,d}:\Sym ^{d,d}(X/S)\to \Sym ^{\infty ,\infty }(X/S)
  $$
be the fibred product of $\iota _d$ with itself over $S$. Recall that $\Sym ^{\infty } (X/S)^+$ is the group completion of the monoid $\Sym ^{\infty }(X/S)$ in the category $\Shv ((\Noe /S)_{\Nis })$. It means that the we have a pushout square
  $$
  \diagram
  \Sym ^{\infty }(X/S)\ar[rr]^-{\Delta }
  \ar[dd]^-{}
  & & \Sym ^{\infty ,\infty }(X/S)
  \ar[dd]^-{\sigma _{\infty }} \\ \\
  S \ar[rr]^-{} & & \Sym (X/S)^+
  \enddiagram
  $$
in the category $\Mon (\Shv ((\Noe /S)_{\Nis }))$. In particular, the quotient morphism $\sigma _{\infty }$ is a morphism of monoids, i.e. it respects the monoidal operations in the source and target. Let
  $$
  \sigma _d:\Sym ^{d,d}(X/S)\to \Sym ^{\infty }(X/S)^+
  $$
be the composition of the morphisms $\iota _{d,d}$ and $\sigma _{\infty }$ in the category $\Shv ((\Noe /S)_{\Nis })$, and let
   $$
   \Sym ^d(X/S)^+
   $$
be the sheaf-theoretical image of the morphism $\sigma _d$, i.e. the image of $\sigma _d$ in the category $\Shv ((\Noe /S)_{\Nis })$.

Some explanation is in place here. A priori, for any nonnegative integer $d$, one can compute the $d$-th symmetric power
  $$
  S^d(X/S)
  $$
in the category of presheaves $\PShv (\Noe /S)$, and the $d$-th symmetric power
  $$
  \Sym ^d(X/S)\; ,
  $$
computed in the category of sheaves $\Shv ((\Noe /S)_{\Nis })$, is the Nisnevich sheafification of the presheaf $S^d(X/S)$. But since the symmetric power $S^d(X/S)$ exists already as a scheme in the category $\Noe /S$, and since the Nisnevich topology is subcanonical, we have that
  $$
  S^d(X/S)=\Sym ^d(X/S)\; ,
  $$
for any $d\geq 0$.

Let
  $$
  \coprod _{d=0}^{\infty }S^d(X/S)
  $$
be the free monoid $\NN (X/S)$ of $X$ over $S$ computed in the category of presheaves $\PShv (\Noe /S)$. Since the category $\Noe /S$ is a Noetherian category, one can show that this infinite coproduct is a Nisnevich sheaf, and hence it coincides with the free monoid $\NN (X/S)$ of $X$ over $S$ computed in the category of sheaves $\Shv ((\Noe /S)_{\Nis })$. In other words, there is no difference between $\NN (X/S)$ in $\PShv (\Noe /S)$ and $\NN (X/S)$ in $\Shv ((\Noe /S)_{\Nis })$, and we write
  $$
  \NN (X/S)=\coprod _{d=0}^{\infty }\Sym ^d(X/S)=
  \coprod _{d=0}^{\infty }S^d(X/S)\; .
  $$

Similarly, let
  $$
  S^{\infty }(X/S)
  $$
be the free connective monoid $\NN (X/S)_{\infty }$ of $X$ over $S$ computed in the category of presheaves $\PShv (\Noe /S)$, so that the free connective monoid $\Sym ^{\infty }(X/S)$ of $X$ over $S$, computed in the category of sheaves $\Shv ((\Noe /S)_{\Nis })$, is nothing but the Nisnevich sheafification of $S^{\infty }(X/S)$. Again, as the category $\Noe /S$ is a Noetherian category, one can show that $S^{\infty }(X/S)$ is a sheaf in Nisnevich topology, and hence
  $$
  S^{\infty }(X/S)=\Sym ^{\infty }(X/S)\; .
  $$

This gives us that, if
  $$
  S^{\infty }(X/S)^+
  $$
is the group completion of the presheaf free monoid $S^{\infty }(X/S)$ in the category $\Mon (\PShv (\Noe /S))$, i.e. the square
  \begin{equation}
  \label{completiondiagr*}
  \diagram
  S^{\infty }(X/S)\ar[rr]^-{\Delta }
  \ar[dd]^-{}
  & & S^{\infty ,\infty }(X/S)
  \ar[dd]^-{\sigma _{\infty }} \\ \\
  S \ar[rr]^-{} & & S(X/S)^+
  \enddiagram
  \end{equation}
is co-Cartesian, the sheaf group completion $\Sym ^{\infty }(X/S)^+$ of $\Sym ^{\infty }(X/S)$ in the category $\Mon (\Shv ((\Noe /S)_{\Nis }))$ is the sheafification of $S^{\infty }(X/S)^+$. In other words,
  $$
  \Sym ^{\infty }(X/S)^+=(S^{\infty }(X/S)^+)^{\shf }\; ,
  $$
where $-^{\shf }$ stays for sheafification functor. 

\begin{lemma}
\label{smallbutimportant}
The presheaf $S^{\infty }(X/S)^+$ is separated. Equivalently, the canonical morphism
  $$
  S^{\infty }(X/S)^+\to \Sym ^{\infty }(X/S)^+
  $$
is a monomorphism in $\PShv (\Noe /S)$.
\end{lemma}

\begin{pf}
Since $S^{\infty }(X/S)^+$ is an abelian group object in the category $\PShv (\Noe /S)$, to prove the lemma it is enough to show that, if
  $$
  F\in S^{\infty }(X/S)^+(U)
  $$
is a section of the presheaf $S^{\infty }(X/S)^+$ on some locally Noetherian scheme $U$ over $S$, and if there exists a Nisnevich covering
  $$
  \{ f_i:U_i\to U\} _{i\in I}\; ,
  $$
such that the pullback $F_i$ of the section $F$ to $U_i$ along each morphism $U_i\to U$ is $0$ in the abelian group $S^{\infty }(X/S)^+(U_i)$, then $F$ is $0$ in the abelian group $S^{\infty }(X/S)^+(U)$.

The section $F$ can be interpreted as a morphism
  $$
  F:U\to S^{\infty }(X/S)^+\; .
  $$
For short of notation, let
  $$
  S^{\infty ,\infty }(X/S)=
  S^{\infty }(X/S)\times _SS^{\infty }(X/S)\; ,
  $$
and, for any nonnegative integer $d$ let
  $$
  S^{d,d}(X/S)=S^d(X/S)\times _SS^d(X/S)\; .
  $$
In these terms, the morphism $F$ is the composition of a certain morphism
  $$
  (f_1,f_2):U\to S^{\infty ,\infty }(X/S)\; ,
  $$
induced by two morphisms of presheaves
  $$
  f_1:U\to S^{\infty }(X/S)
  \qqand
  f_1:U\to S^{\infty }(X/S)\; ,
  $$
and the quotient morphism
  $$
  \sigma _{\infty }:S^{\infty ,\infty }(X/S)\to
  S^{\infty }(X/S)^+\; .
  $$
Moreover, there exists $d$, such that both morphisms $f_1$ and $f_2$ factorize through $S^d(X/S)$, and then $F$ is the composition
  \begin{equation}
  \label{zimodry}
  U\stackrel{(f_1,f_2)}{\lra }S^{d,d}(X/S)
  \stackrel{\iota _{d,d}}{\lra }S^{\infty ,\infty }(X/S)
  \stackrel{\sigma _{\infty }}{\lra }S^{\infty }(X/S)^+\; .
  \end{equation}
and the morphisms $f_1$ and $f_2$ are morphisms of locally Noetherian schemes over the base scheme $S$.

Now, since $S^{\infty }(X/S)$ is a cancellative monoid in $\PShv (\Noe /S)$, the commutative square (\ref{completiondiagr*}) is a Cartesian square in $\PShv (\Noe /S)$. It follows that, since $F_i=0$ for all $i\in I$, the images of the compositions
  $$
  U_i\to U\stackrel{(f_1,f_2)}{\lra }
  S^{d,d}(X/S)\stackrel{\iota _{d,d}}{\lra }
  S^{\infty ,\infty }(X/S)
  \stackrel{\sigma _{\infty }}{\lra }S^{\infty }(X/S)^+
  $$
are all in the image of the diagonal morphism
  $$
  \Delta :S^{\infty }(X/S)\to S^{\infty ,\infty }(X/S)\; .
  $$
And since the morphism
  $$
  \coprod _{i\in I}U_i\to U
  $$
is a scheme-theoretical epimorphism, we see that the image of the morphism (\ref{zimodry}) is also in the image of the diagonal morphism $\Delta $. The latter means that the section $F$ equals $0$.
\end{pf}

Let
  $$
  \sigma _d:S^{d,d}(X/S)\to S^{\infty }(X/S)^+
  $$
be the composition of the morphisms $\iota _{d,d}$ and $\sigma _{\infty }$ in the category $\PShv (\Noe /S)$, and let
   $$
   S^d(X/S)^+
   $$
be the image of the morphism $\sigma _d$ in the category $\PShv (\Noe /S)$. Then $S^d(X/S)^+$ is a sub-presheaf in $\Sym ^{\infty }(X/S)^+$. As the sheafification functor is exact, it preserves monomorphisms. It follows that
   $$
   \Sym ^d(X/S)^+=(S^d(X/S)^+)^{\shf }\; ,
   $$
i.e. $\Sym ^d(X/S)^+$ is the Nisnevich sheafification of the preshaef $S^d(X/S)^+$. And, once again, the sheaf-theoretical image $\Sym ^d(X/S)^+$ of the morphism $\sigma _d$ comes together with the epimorphism
  \begin{equation}
  \label{vechervgrumbi_z}
  \sigma _d:\Sym ^{d,d}(X/S)\to \Sym ^d(X/S)^+
  \end{equation}
in the category $\Shv ((\Noe /S)_{\Nis })$.

Next, the section $S\to X$ of the structural morphism $X\to S$ induces the closed embeddings
  $$
  \Sym ^d(X/S)\to \Sym ^{d+1}(X/S)\; ,
  $$
which, in turn, induce the closed embeddings
  $$
  \Sym ^{d,d}(X/S)\to \Sym ^{d+1,d+1}(X/S)\; .
  $$
The latter morphisms induce the corresponding morphisms
  $$
  \Sym ^d(X/S)^+\to \Sym ^{d+1}(X/S)^+
  $$
in the category $\Shv ((\Noe /S)_{\Nis })$. Then
  \begin{equation}
  \label{indstructure1}
  \Sym ^d(X/S)^+=\colim _d\; \Sym ^d(X/S)^+\; ,
  \end{equation}
i.e. the space $\Sym ^d(X/S)^+$ is naturally the colimit of the spaces $\Sym ^d(X/S)^+$.

\begin{remark}
{\rm The sheaf $\Sym ^d(X/S)^+$ is {\it not} a group completion of any monoid.
}
\end{remark}

The constructions above allow us to consider a natural atlas for the
  $$
  CA_0(X/S,0)=\{ \sigma _d\; |\; d\in \ZZ \; ,\; d\geq 0\}
  $$
be the set of all morphisms $\sigma _d$, and let
  $$
  CA(X/S,0)=\langle CA_0(X/S,0)\rangle
  $$
be the {\it Chow atlas} on the Nisnevich connective space $\Sym ^{\infty }(X/S)^+$. According to Section \ref{kaehler}, the sheaf $\Sym ^{\infty }(X/S)^+$ is now the Nisnevich space of relative $0$-cycles on $X$ over $S$, with regard to the Chow atlas
  $$
  CA=CA(X/S,0)\; .
  $$
For short, we will say that $\Sym ^{\infty }(X/S)^+$ is the {\it space of $0$-cycles} on $X$ over $S$

Hilbert schemes allow us to consider a natural subatlas in the Chow atlas $CA$. Indeed, let $U$ be a locally Noetherian scheme over $S$, and let
  $$
  Z\to X\times _SU
  $$
be a closed subscheme in $X\times _SU$. Suppose the composition
  $$
  g:Z\to U
  $$
of the closed embedding of $Z$ into $X\times _SU$ with the projection onto $U$ is flat. Then, if $V$ is an irreducible component of $Z$, the closure $\overline {g(V)}$ is an irreducible component of $U$. Therefore, if $U$ is irreducible, $\overline {g(V)}=U$. If, moreover, $g$ is proper, then $\overline {g(V)}=g(V)$, and hence $g$ is a surjection.

Since $X$ is embedded in to $\PR ^n_S$ over $S$ via the closed embedding $i$, the scheme $X\times _SU$ embeds into $\PR ^m_U$ over $U$, and the morphism $g:Z\to U$ factorizes through the embedding of $Z$ into $\PR ^m_U$ followed by the projection from $\PR ^m_U$ onto $U$. Therefore, if $u\in U$ and $Z_u$ is the fibre of $g$ at $u$, the Hilbert polynomial of the structural sheaf $\bcO _{Z_u}$ does not depend on $u$, see Theorem 9.9 on page 261 in \cite{Hartshorne}.

This fact allows us to consider, for every polynomial
  $$
  P\in \QQ [x]
  $$
the standard Hilbert set valued presheaf
  $$
  \HilbF _P(X/S):\Noe /S\to \Sets
  $$
sending a locally Noetherian $S$-scheme $U$ to the set of closed subschemes $Z$ in the product $X\times _SU$, which are flat and proper over $U$, and such that the Hilbert polynomial of $\bcO _{Z_u}$ is $P$. Let also
  $$
  \HilbF(X/S)=\coprod _{P\in \QQ [x]}\HilbF _P(X/S):
  \Noe /S\to \Sets
  $$
be the total Hilbert functor on locally Noetherian schemes over $S$.

Since $X$ is projective over $S$, the Hilbert functors $\HilbF _P(X/S)$ are representable. This result is due to Grothendieck, see Chapter 5 in \cite{FGAexplained} or Chapter I.1 in \cite{KollarRatCurvesOnVar}. For each polynomial $P$ in $\QQ [x]$ there exists a scheme, called the Hilbert scheme,
  $$
  \HilbS _P(X/S)
  $$
over $S$ representing the functor $\HilbF _P(X/S)$. Moreover, this scheme is projective over $S$.

Within this paper we are interested in the case when $P=d$ is a non-negative integer. In that case the Hilbert scheme
  $$
  \HilbS ^d(X/S)=\HilbS _P(X/S)|_{P=d}
  $$
is a scheme over the $d$-th relative symmetric power, and we have the so-called Hilbert-Chow morphism of schemes
   \begin{equation}
   \label{enot_d}
   \hc _d:\HilbS ^d(X/S)\to \Sym ^d(X/S)\; .
   \end{equation}

For any nonnegative integer $d$ let
  $$
  \HilbS ^{d,d}(X/S)=
  \HilbS ^d(X/S)\times _S\HilbS ^d(X/S)\; ,
  $$
and let
  $$
  HA_0(X/S,0)=\{ a_d\circ (\hc _{d,d})\; |\;
  d\in \ZZ \; ,\; d\geq 0\} \; ,
  $$
where
  $$
  \hc _{d,d}:\HilbS ^{d,d}(X/S)\to
  \Sym ^{d,d}(X/S)
  $$
is the fibred self-product over $S$ of the $d$-th Hilbert-Chow morphism $\hc _d$. Let also
  $$
  HA(X/S,0)=\langle HA_0(X/S,0)\rangle
  $$
be the {\it Hilbert atlas} on the space $\Sym ^{\infty }(X/S)^+$. Obviously, the Hilbert atlas is a subatlas of the Chow atlas on $\Sym ^{\infty }(X/S)^+$.

Now, let
  $$
  \bcO _{\Sym ^{\infty }(X/S)^+}
  $$
be the sheaf of regular functions on the site $\Sym ^{\infty }(X/S)^+_{\Nis \mhyphen \et }$, constructed with regard to the Chow atlas $CA$ on the sheaf $\Sym ^{\infty }(X/S)^+$, as explained in Section \ref{kaehler}. In particular, if $U\to \Sym ^{\infty }(X/S)^+$ is a morphism from a scheme $U$ to $\Sym ^{\infty }(X/S)^+$ over $S$, which is \'etale with regard to the Chow atlas on $\Sym ^{\infty }(X/S)^+$, then
since
  $$
  \Gamma (U\to \Sym ^{\infty }(X/S)^+,
  \bcO _{\Sym ^{\infty }(X/S)^+})=
  \Gamma (U,\bcO _U)\; .
  $$

As soon as the sheaf $\bcO _{\Sym ^{\infty }(X/S)^+}$ is defined, we can also define the sheaf of K\"ahler differentials
  $$
  \Omega ^1_{\Sym ^{\infty }(X/S)^+}=
  \Omega ^1_{\Sym ^{\infty }(X/S)^+/S}
  $$
on the site $\Sym ^{\infty }(X/S)^+_{\Nis \mhyphen \et }$, see Section \ref{kaehler}. Let also
  $$
  T_{\Sym ^{\infty }(X/S)^+}=
  T_{\Sym ^{\infty }(X/S)^+/S}
  $$
be the tangent sheaf, i.e. the dual to the sheaf of K\"ahler differentials on the site $\Sym ^{\infty }(X/S)^+_{\Nis \mhyphen \et }$.

Since now the sheaf of K\"ahler differentials and the tangent sheaf on the site $\Sym ^{\infty }(X/S)^+_{\Nis \mhyphen \et }$ will be considered as the sheaf of K\"ahler differentials and the tangent sheaf on the space of $0$-cycles $\Sym ^{\infty }(X/S)^+$.

Notice that both sheaves are given in terms of the Chow atlas on $\Sym ^{\infty }(X/S)^+$. Similar sheaves can be also defined in terms of the Hilbert atlas on the same space, and the connection between two types is an interesting question, also considered in \cite{GreenGriffiths}, but in different terms.

\subsection{\'Etale neigbourhoods of a point on $\Sym ^{\infty }(X/S)^+$}

Next, recall that a point $P$ on $\Sym ^{\infty }(X/S)^+$ is an equivalence class of morphisms from spectra of fields to $\Sym ^{\infty }(X/S)^+$, as explained in Section \ref{kaehler}. By abuse of notation, we write
  $$
  P:\Spec (K)\to \Sym ^{\infty }(X/S)^+\; .
  $$
We will always assume that $P$ factorizes through the Chow atlas $CA$ on the space $\Sym ^{\infty }(X/S)^+$.

As in Section \ref{kaehler}, consider the functor
  $$
  u_P:\Sym ^{\infty }(X/S)^+_{\Nis \mhyphen \et }\to \Sets
  $$
sending an \'etale morphism
  $$
  U\to \Sym ^{\infty }(X/S)^+\; ,
  $$
where $U$ is a locally Noetherian scheme over $S$, to the set
  $$
  u_P(U)=|U_P|
  $$
of points on the fibre
  $$
  U_P=U\times _{\Sym ^{\infty }(X/S)^+}\Spec (K)
  $$
of the morphism $U\to \Sym ^{\infty }(X/S)^+$ at $P$.

As soon as the functor $u_P$ is introduced, we also define the notion of a neighbourhood of $P$, with regard to the functor $u_P$, as we did it in Section \ref{kaehler}. Namely, an \'etale neighbourhood of $P$ on $\Sym ^{\infty }(X/S)^+$ is a pair
  $$
  N=(U\to \Sym ^{\infty }(X/S)^+,T\in u_P=|U_P|)\; ,
  $$
where the morphism $U\to \Sym ^{\infty }(X/S)^+$ is over $S$ and \'etale with regard to the Chow atlas $CA$ on $\Sym ^{\infty }(X/S)^+$, and $T$ is a point of the fibre $U_P$. Or, equivalently, an \'etale neighbourhood of $P$ is an \'etale morphism
  $$
  U\to \Sym ^{\infty }(X/S)^+
  $$
over $S$ such that the point
  $$
  P:\Spec (K)\to \Sym ^{\infty }(X/S)^+
  $$
factorizes through $U$.

As in Section \ref{kaehler}, all \'etale neighbourhoods form the category of \'etale neighbourhoods of $P$ on $\Sym ^{\infty }(X/S)^+$ denoted by $\bcN _P$.

Now, Tag 00YA \cite{StacksProject} gives us that in order to show that the corresponding stalk functor
  $$
  \stalk _P:
  \Shv (\Sym ^{\infty }(X/S)^+_{\Nis \mhyphen \et })\to
  \Sets
  $$
induces a point of the topos $\Shv (\Sym ^{\infty }(X/S)^+_{\Nis \mhyphen \et })$, we need to show that the functor $u_P$ satisfies all the three items of Tag 00Y5 in \cite{StacksProject}. The items (1) and (2) are satisfied in general, see Section \ref{kaehler}. The last item (3) of Tag 00Y5 in \cite{StacksProject} is satisfied when the category $\bcN _P$ is cofiltered. Therefore, our aim is now to show that, in case of the space of $0$-cycles $\Sym ^{\infty }(X/S)^+$ the category $\bcN _P$ is cofiltered.

We start with the following representability lemma, which will be necessary for the study of the category $\bcN _P$.

\begin{lemma}
\label{keylemma}
For any nonnegative integer $d$ and for any two morphisms
  $$
  U\to S^{\infty }(X/S)^+
  \qqand
  V\to S^{\infty }(X/S)^+\; ,
  $$
where $U$ and $V$ are locally Noetherian schemes over $S$, the fibred product
  $$
  U\times _{S^{\infty }(F/S)^+}V\; ,
  $$
in the category of presheaves $\PShv (X/S)$, is represented by a locally Noetherian scheme over $S$.
\end{lemma}

\begin{pf}
We need to find a locally Noetherian scheme over $S$ representing the fibred product
  $$
  U\times _{S^{\infty }(X/S)^+}V
  $$
in the category $\PShv (\Noe /S)$.

Denote the morphism from $U$ to $S^{\infty }(X/S)^+$ by $F$, and the morphism from $V$ to $S^{\infty }(X/S)^+$ by $G$.
\iffalse %%%%%%%%%%%%%%%%%%%%%%%%%%%%%%%%%%%%%%%%%%%%%%%%%
  $$
  \diagram
  U\times _SV \ar@<0.5ex>[dddrrr]^-{G\circ \pr _V}
  \ar@<-0.5ex>[dddrrr]_-{F\circ \pr _U}
  \ar[dddd]_-{\pr _U} \ar[rrrr]^-{\pr _V} & & & &
  V \ar[dddl]_-{G}\ar[dddd]^-{} \\ \\ \\
  & & & S^{\infty }(X/S)^+ \ar[rd]^-{} & \\
  U \ar[rrru]^-{F} \ar[rrrr]^-{}
  & & & & S
  \enddiagram
  $$
\fi %%%%%%%%%%%%%%%%%%%%%%%%%%%%%%%%%%%%%%%%%%%%%%%%%%%%%%
Clearly, the object $U\times _{S^{\infty }(F/S)^+}V$ is the coequalizer of the compositions of the projections from $U\times _SV$ on to $U$ and $V$ with the morphisms $F$ and $G$ respectively.
\iffalse %%%%%%%%%%%%%%%%%%%%%%%%%%%%%%%%%%%%%%%%%%%%%%%%%
  $$
  U\times _{\Sym ^{\infty }(F/S)^+}V=
  \coeq (\xymatrix{
  U\times _SV \ar@<+0.5ex>[r]^-{} \ar@<-0.5ex>[r]^-{} &
  S^{\infty }(X/S)^+})
  $$
\fi %%%%%%%%%%%%%%%%%%%%%%%%%%%%%%%%%%%%%%%%%%%%%%%%%%%%%%
Since $S^{\infty }(X/S)^+$ is an abelian group object, one can consider the difference
  $$
  H=F\circ \pr _U-G\circ \pr _V:U\times _SV\to
  S^{\infty }(X/S)^+
  $$
between these two compositions in the category $\PShv (\Noe /S)$. Then the coequalizer $U\times _{\Sym ^{\infty }(F/S)^+}V$ fits in to the Cartesian square
  $$
  \diagram
  U\times _{\Sym ^{\infty }(F/S)^+}V
  \ar[dd]_-{} \ar[rr]^-{} & & U\times _SV
  \ar[dd]^-{H} \\ \\
  S \ar[rr]^-{} & & S^{\infty }(X/S)^+
  \enddiagram
  $$
and the lemma reduces to the case when $U=S$, and $F$ is a section of the structural morphism from $S^{\infty }(X/S)^+$ to $S$.

Next, the morphism of presheaves
  $$
  G:V\to S^{\infty }(X/S)^+
  $$
is uniquely determined by sending the identity morphism $\id _V$ to some element in the abelian group $S^{\infty }(X/S)^+(V)$, which is the equivalence class
   $$
   [(g_1,g_2)]
   $$
of two morphisms
   $$
   g_1:V\to S^{\infty }(X/S)
   \qqand
   g_2:V\to S^{\infty }(X/S)
   $$
of presheaves over $S$. In particular, the morphism $G$ factorized through the product
  $$
  S^{\infty ,\infty }(X/S)=
  S^{\infty }(X/S)\times _SS^{\infty }(X/S)\; .
  $$
As we mentioned already, since $S^{\infty }(X/S)$ is a cancellation monoid, the commutative square (\ref{completiondiagr*}) is a Cartesian square in $\PShv (\Noe /S)$. Let
  $$
  V_S=S^{\infty }(X/S)\times _{S^{\infty ,\infty }(X/S)}V
  $$
be the fibred product of $S^{\infty }(X/S)$ and $V$ over $S^{\infty ,\infty }(X/S)$. The composition of the two Cartesian squares
    $$
    \xymatrix{
    V_S \ar[dd]_-{} \ar[rr]^-{} & & V \ar[dd]^-{} \\ \\
    S^{\infty }(X/S) \ar[rr]^-{\Delta } & &
    S^{\infty ,\infty }(X/S)
    }
    $$
and
    $$
    \xymatrix{
    S^{\infty }(X/S) \ar[rr]^-{}
    \ar[dd]^-{} & & S^{\infty ,\infty }(X/S)
    \ar[dd]^-{} \\ \\
    S \ar[rr]^-{} & & S^{\infty }(X/S)^+
    }
    $$
shows that the object $V_S$ fits in to the Cartesian square
    $$
    \xymatrix{
    V_S \ar[dd]_-{} \ar[rr]^-{} & & V \ar[dd]^-{G} \\ \\
    S \ar[rr]^-{} & & S^{\infty }(X/S)^+
    }
    $$
In other words,
  $$
  V_S=S\times _{S^{\infty }(X/S)^+}V
  $$
is, at the same time, the fibred product of $S$ and $V$ over $S^{\infty }(X/S)^+$.

Choose $d$ such that the image of the morphism
  $$
  V\to S^{\infty ,\infty }(X/S)
  $$
is in
  $$
  S^{d,d}(X/S)=S^d(X/S)\times _SS^d(X/S)\; .
  $$
Since the morphism
  $$
  S^{d,d}(X/S)\to S^{\infty ,\infty }(X/S)
  $$
is a monomorphism in $\PShv (\Noe /S)$, it follows that the object $V_S$ fits also in to the Cartesian square
    $$
    \xymatrix{
    V_S \ar[dd]_-{} \ar[rr]^-{} & & V \ar[dd]^-{} \\ \\
    S^d(X/S) \ar[rr]^-{} & & S^{d,d}(X/S)
    }
    $$
In other words, $V_S$ is the fibred product of the schemes $S^d(X/S)$ and $V$ over the scheme $S^{d,d}(X/S)$. In particular, $V_S$ is a scheme itself.
\end{pf}

We need one easy but useful technical notion. Suppose we are given with a locally Noetherian scheme $U$ over $S$ and a morphism
   $$
   F:U\to \Sym ^{\infty }(X/S)^+
   $$
in the category of sheaves $\Shv ((\Noe /S)_{\Nis })$. Any such a morphism is uniquely determined by sending $\id _U:U\to U$ to a section
  $$
  s_F\in \Sym ^{\infty }(X/S)^+(U)\; .
  $$
Since $\Sym ^{\infty }(X/S)^+$ is the Nisnevich sheafification of the presheaf $S^{\infty }(X/S)^+$ and the morphism
  $$
  S^{\infty }(X/S)^+\to \Sym ^{\infty }(X/S)^+
  $$
is a monomorphism in $\PShv (\Noe /S)$ by Lemma \ref{smallbutimportant}, we obtain that the section $s_F$ is the equivalence class of pairs, each of which consists of a Nisnevich cover
  $$
  \{ U_i\to U\} _{i\in I}
  $$
and a collection of sections
  $$
  s_i\in S^{\infty }(X/S)^+(U_i)\; ,
  $$
such that the restrictions of $s_i$ and $s_j$ on $U_i\times _UU_j$ coincide for all indices $i$ and $j$ in $I$. Therefore, if $$
  \hat U=\coprod _{i\in I}U_i\; ,
  $$
we obtain two morphisms
  $$
  \hat U\to U
  $$
and
  $$
  \hat F:\hat U\to S^{\infty }(X/S)^+
  $$
such that the square
  $$
  \diagram
  \hat U\ar[dd]_-{} \ar[rr]^-{\hat F} & &
  S^{\infty }(X/S)^+ \ar[dd]^-{} \\ \\
  U \ar[rr]^-{F} & & \Sym ^{\infty }(X/S)^+
  \enddiagram
  $$
is commutative. For short, we will say that $\hat U$ (respectively, $\hat F$) is an extension of $U$ (resp., $F$) by the pair $(\{ U_i\to U\} _{i\in I},\{ s_i\} _{i\in I})$ representing the section $s_F$.

\subsection{The category $\bcN _P$ and the proof of Main Theorem}
\label{mainthm}

Now we are ready to prove our main result.

\begin{theorem}
\label{cofilter}
Let $P$ be a point of the space $\Sym ^{\infty }(X/S)^+$, and let $\bcN _P$ be the category of \'etale neighbourhoods of the point $P$ on $\Sym ^{\infty }(X/S)^+$. Then $\bcN _P$ is cofiltered.
\end{theorem}

\begin{pf}
The proof follows a pretty standard way of argumentation, see, for example, Lemma 57.18.3. First of all, Lemma \ref{keylemma} gives us that the category $\bcN _P$ is nonempty, so that the first axiom of a cofiltered category is satisfied.

Let
  $$
  F:U\to \Sym ^{\infty }(X/S)^+
  \qqand
  G:V\to \Sym ^{\infty }(X/S)^+
  $$
\iffalse %%%%%%%%%%%%%%%%%%%%%%%%%%%%%%%%%%%%%%%%%%%%%%%%%%%%
    $$
    \xymatrix{
    \Spec (K) \ar[dd]^-{} \ar[ddrr]^-{P} & & &
    \Spec (K) \ar[dd]^-{} \ar[ddrr]^-{P} & & \\ \\
    U \ar[rr]^-{F} & & \Sym ^{\infty }(X/S)^+ &
    V \ar[rr]^-{G} & & \Sym ^{\infty }(X/S)^+}
    $$
\fi %%%%%%%%%%%%%%%%%%%%%%%%%%%%%%%%%%%%%%%%%%%%%%%%%%%%%%%%%
be two \'etale neighbourhoods of the point $P$, and look at the fibred product
  $$
  \diagram
  U\times _{\Sym ^{\infty }(X/S)^+}V
  \ar[dd]_-{} \ar[rr]^-{} & & V \ar[dd]^-{G} \\ \\
  U \ar[rr]^-{F} & & \Sym ^{\infty }(X/S)^+
  \enddiagram
  $$

Let $s_F$ be the section of the sheaf $\Sym ^{\infty }(X/S)^+$ on $U$ which determines the morphism $F$, and let
  $$
  \hat F:\hat U\to S^{\infty }(X/S)^+
  $$
be the extension of the morphism $F$ given by a pair $(\{ U_i\to U\} _{i\in I},\{ s_i\} _{i\in I})$ representing $s_F$. Similarly, one can construct an extension $\hat G$ of the morphism $G$ induced by a pair representing the section $s_G$.

By Lemma \ref{keylemma}, the fibred product $\hat U\times _{S^{\infty }(X/S)^+}\hat V$ is a locally Noetherian scheme over $S$. Consider the universal morphism
  $$
  \hat U\times _{S^{\infty }(X/S)^+}\hat V\to
  U\times _{\Sym ^{\infty }(X/S)^+}V\; ,
  $$
commuting with the extensions of $U$ and $V$.
\iffalse %%%%%%%%%%%%%%%%%%%%%%%%%%%%%%%%%%%%%%%%%%%%%%%%%%%%
  $$
  \diagram
  \hat U\times _{S^{\infty }(X/S)^+}\hat V \ar[dddd]^-{}
  \ar[rr]^-{}
  \ar[ddr]^-{} & &
  \hat V \ar[dddd]^-{}
  \ar[ddr]^-{} & & \\ \\
  & U\times _{\Sym ^{\infty }(X/S)^+}V \ar[dddd]^-{}
  \ar[rr]^-{}
  & & V \ar[dddd]^-{G} & \\ \\
  \hat U \ar[rr]^-{}
  \ar[ddr]^-{} & & S^{\infty }(X/S)^+
  \ar[ddr]^-{} & & \\ \\
  & U \ar[rr]^-{F}
  & & \Sym ^{\infty }(X/S)^+ &
  \enddiagram
  $$
\fi %%%%%%%%%%%%%%%%%%%%%%%%%%%%%%%%%%%%%%%%%%%%%%%%%%%%%%%%%
Let us show that the composition
  $$
  H:\hat U\times _{S^{\infty }(X/S)^+}\hat V\to
  U\times _{\Sym ^{\infty }(X/S)^+}V\to
  \Sym ^{\infty }(X/S)^+
  $$
is \'etale, with regard to the Chow atlas on $\Sym ^{\infty }(X/S)^+$.

Indeed, since the morphism
  $$
  S^{\infty }(X/S)^+\to \Sym ^{\infty }(X/S)^+
  $$
is a monomorphism in $\PShv (\Noe /S)$ by Lemma \ref{smallbutimportant}, the square
  $$
  \diagram
  \hat U\times _{S^{\infty }(X/S)^+}\hat V
  \ar[dd]_-{\id } \ar[rr]^-{} & & S^{\infty }(X/S)^+
  \ar[dd]^-{} \\ \\
  \hat U\times _{S^{\infty }(X/S)^+}\hat V
  \ar[rr]^-{} & & \Sym ^{\infty }(X/S)^+
  \enddiagram
  $$
is Cartesian, so that the obvious morphism
  $$
  h:\hat U\times _{S^{\infty }(X/S)^+}\hat V\to
  S^{\infty }(X/S)^+
  $$
is the pullback of the morphism $H$.

For short, let
  $$
  \hat U_{d,d}=
  \hat U\times _{S^{\infty }(X/S)^+}S^{d,d}(X/S)\; ,
  $$
and
  $$
  \hat V_{d,d}=
  \hat V\times _{S^{\infty }(X/S)^+}S^{d,d}(X/S)\; .
  $$
Then
  $$
  h_0:\hat U_{d,d}\times _{S^{d,d}(X/S)}\hat V_{d,d}
  \to S^{d,d}(X/S)
  $$
is the pullback of the morphism $h$, and since $h$ is the pullback of $H$, we obtain the Cartesian square
  $$
  \diagram
  \hat U_{d,d}\times _{S^{d,d}(X/S)}\hat V_{d,d}
  \ar[dd]_-{} \ar[rr]^-{h_0} & & S^{d,d}(X/S) \ar[dd]^-{} \\ \\
  \hat U\times _{S^{\infty }(X/S)^+}\hat V
  \ar[rr]^-{H} & & \Sym ^{\infty }(X/S)^+
  \enddiagram
  $$
\iffalse %%%%%%%%%%%%%%%%%%%%%%%%%%%%%%%%%%%%%%%%%%%%%%%%%%%%%
  $$
  \diagram
  \hat U_{d,d}\times _{S^{d,d}(X/S)}\hat V_{d,d}
  \ar[dddd]^-{} \ar[rr]^-{} \ar[ddr]^-{} & &
  \hat V_{d,d} \ar[dddd]^-{} \ar[ddr]^-{} & & \\ \\
  & \hat U\times _{S^{\infty }(X/S)^+}\hat V \ar[dddd]^-{}
  \ar[rr]^-{}
  & & \hat V \ar[dddd]^-{\hat G} & \\ \\
  \hat U_{d,d} \ar[rr]^-{}
  \ar[ddr]^-{} & & S^{d,d}(X/S)
  \ar[ddr]^-{} & & \\ \\
  & \hat U \ar[rr]^-{\hat F}
  & & S^{\infty }(X/S)^+ &
  \enddiagram
  $$
\fi %%%%%%%%%%%%%%%%%%%%%%%%%%%%%%%%%%%%%%%%%%%%%%%%%%%%%%%%%%%
Therefore, in order to prove that $H$ is \'etale, we need only to show that $h_0$ is \'etale.

Now again, since the morphism from $S^{\infty }(X/S)^+$ to $\Sym ^{\infty }(X/S)^+$ is a monomorphism in $\PShv (\Noe /S)$ by Lemma \ref{smallbutimportant}, we see that the commutative square
  $$
  \diagram
  \hat U \ar[dd]_-{\id } \ar[rr]^-{} & & S^{\infty }(X/S)^+ \ar[dd]^-{} \\ \\
  \hat U \ar[rr]^-{} & & \Sym ^{\infty }(X/S)^+
  \enddiagram
  $$
is Cartesian. Composing it with the Cartesian square
  $$
  \diagram
  \hat U_{d,d} \ar[dd]_-{} \ar[rr]^-{} & & S^{d,d}(X/S) \ar[dd]^-{} \\ \\
  \hat U \ar[rr]^-{} & & S^{\infty }(X/S)^+
  \enddiagram
  $$
we obtain the Cartesian square
  $$
  \diagram
  \hat U_{d,d} \ar[dd]_-{} \ar[rr]^-{} & & S^{d,d}(X/S) \ar[dd]^-{} \\ \\
  \hat U \ar[rr]^-{} & & \Sym ^{\infty }(X/S)^+
  \enddiagram
  $$

The bottom horizontal morphism is the composition of two \'etale morphisms, and hence it is \'etale. Since \'etale morphisms are stable under pullbacks, the top horizontal morphism
  $$
  \hat U_{d,d}\to S^{d,d}(X/S)
  $$
in the latter square is \'etale as well.

Similarly, the morphism
  $$
  \hat V_{d,d}\to S^{d,d}(X/S)
  $$
is \'etale.

Thus, the bottom horizontal and the right vertical morphisms in then Cartesian square
  $$
  \diagram
  \hat U_{d,d}\times _{S^{d,d}(X/S)}\hat V_{d,d}
  \ar[dd]_-{} \ar[rr]^-{} & & \hat V_{d,d} \ar[dd]^-{} \\ \\
  \hat U_{d,d} \ar[rr]^-{} & & S^{d,d}(X/S)^+
  \enddiagram
  $$
are \'etale. Since \'etale morphisms are stable under pullbacks and compositions, the diagonal composition $h_0$ of this square is \'etale as well.

As this is true for any $d$, we see that the morphism
  $$
  \hat U\times _{S^{\infty }(X/S)^+}\hat V\to
  \Sym ^{\infty }(X/S)^+
  $$
is \'etale.

The fact that the point $P:\Spec (K)\to \Sym ^{\infty }(X/S)^+$ factorizes through $\hat U\times _{S^{\infty }(X/S)^+}\hat V$ is obvious.

Now we need to prove the last axiom of a cofiltered category. Assume again that we have two \'etale neighbourhoods $U$ and $V$ of $P$ as above, and assume also that we have two morphisms
  $$
  \xymatrix{a,b:U \ar@<+0.7ex>[r]^-{} \ar@<-0.1ex>[r]^-{} & V}
  $$
\iffalse %%%%%%%%%%%%%%%%%%%%%%%%%%%%%%%%%%%%%%%%%%%%%%%%%%%%%%
   \begin{equation}
   \label{reppa1}
   \xymatrix{
   & & \Spec (K) \ar[lldd]_-{} \ar[rrdd]^-{} & & \\ \\
   U \ar@<+0.5ex>[rrrr]^-{a} \ar@<-0.5ex>[rrrr]_-{b}
   \ar[rrdd]^-{F} & & & &
   V  \ar[lldd]_-{G} \\ \\
   & & \Sym ^{\infty }(X/S)^+ & &
   }
   \end{equation}
\fi %%%%%%%%%%%%%%%%%%%%%%%%%%%%%%%%%%%%%%%%%%%%%%%%%%%%%%%%%%%
between these neighbourhoods.

Let
  $$
  s_G\in \Sym ^{\infty }(X/S)^+(V)
  $$
be the section determined by the morphism $G$, and choose a representative in $s_G$. Such a representative consists of a Nisnevich covering
  $$
  \{ V_i\to V\} _{i\in I}
  $$
and a collection of sections
  $$
  s_i\in S^{\infty }(X/S)^+(V_i)\; ,
  $$
such that the restrictions of $s_i$ and $s_j$ on $V_i\times _VV_j$ coincide for all indices $i$ and $j$ in $I$. Construct the corresponding extension
  $$
  \hat G:\hat V\to S^{\infty }(X/S)^+
  $$
of the morphism $G$ getting the commutative square
  $$
  \diagram
  \hat V \ar[dd]_-{} \ar[rr]^-{\hat G} & &
  S^{\infty }(X/S) \ar[dd]^-{} \\ \\
  V \ar[rr]^-{G} & & \Sym ^{\infty }(X/S)^+
  \enddiagram
  $$
Pulling back the \'etale covering $\{ V_i\to V\} _{i\in I}$ along the morphisms $a$ and $b$, and taking the unification
  $$
  \{ U_{ij}\to U\} _{(i,j)\in I\times I}
  $$
of these two pullback coverings in to one, one can construct the extension
  $$
  \hat F:\hat U\to S^{\infty }(X/S)^+\; ,
  $$
such that the diagram
  \begin{equation}
  \label{reppa2}
  \diagram
  \Spec (K)^{\hat {\; }} \ar[dd]_-{} \ar[rr]_-{}
  & & \hat V \ar[dd]^-{\hat G} \\ \\
  \hat U \ar@<+0.5ex>[rruu]^-{\hat a}
  \ar@<-0.5ex>[rruu]_-{\hat b}
  \ar[rr]^-{\hat F} & & S^{\infty }(X/S)^+
  \enddiagram
  \end{equation}
is commutative, where $\Spec (K)^{\hat {\; }}$ is an extension over $\Spec (K)$. Moreover, the squares
    $$
    \xymatrix{
    \hat U \ar[rr]^-{\hat a} \ar[dd]^-{} & &
    \hat V \ar[dd]^-{} & &
    \hat U \ar[dd]^-{} \ar[rr]^-{\hat b} & &
    \hat V \ar[dd]^-{} \\ \\
    U \ar[rr]^-{a} & & V & &
    U \ar[rr]^-{b} & & V}
    $$
    $$
    \xymatrix{
    \Spec (K)^{\hat {\; }} \ar[rr]^-{} \ar[dd]^-{} & &
    \Spec (K) \ar[dd]^-{} &
    \Spec (K)^{\hat {\; }} \ar[dd]^-{} \ar[rr]^-{} & &
    \Spec (K) \ar[dd]^-{} \\ \\
    \hat U \ar[rr]^-{} & & U &
    \hat V \ar[rr]^-{} & & V}
    $$
are commutative.

Now, let $W$ be the fibred product of $U$ and $V$ over $V\times _{\Sym ^{\infty }(X/S)^+}V$, with regard to the morphisms $(a,b)$ and $\Delta $, and let $h$ be the corresponding universal morphism, as it is shown in the commutative diagram
    \begin{equation}
    \label{muhomory1a}
    \xymatrix{
    \Spec (K)\ar@/_/[dddr]^-{} \ar@/^/[drrr]^-{}
    \ar@{.>}[dr]^-{\hspace{-1mm}\exists ! h} \\
    & W \ar[dd] \ar[rr] & & V \ar[dd]^-{\Delta } \\ \\
    & U \ar[rr]^-{(a,b)} & & V\times _{\Sym ^{\infty }(X/S)^+}V}
    \end{equation}
Notice that the external commutativity is guaranteed by the fact that $a$ and $b$ are two morphisms from the neigbourhood $U$ to the neigbourhood $V$ of the same point $P$. The diagram (\ref{muhomory1a}) can be also extended by the commutative diagram
  \begin{equation}
  \label{muhomory1b}
  \diagram
  V\times _{\Sym ^{\infty }(X/S)^+}V
  \ar[dd]^-{} \ar[rr]^-{} & & V \ar[dd]^-{G} \\ \\
  V \ar[rr]^-{G} & & \Sym ^{\infty }(X/S)^+
  \enddiagram
  \end{equation}

\iffalse %%%%%%%%%%%%%%%%%%%%%%%%%%%%%%%%%%%%%%%%%%%%%%%%%%%%%%
Summarizing, we obtain the commutative diagram
  \begin{equation}
  \label{muhomory1}
  \diagram
  \Spec (K)\ar[rrrd]^-{}
  \ar@{.>}[rd]_-{\hspace{+10mm}\exists ! h}
  \ar[rddd]^-{} & & & & & \\
  & W \ar[dd]_-{} \ar[rr]^-{} & & V \ar[rrdd]^-{\id }
  \ar[dd]^-{\Delta } & & \\ \\
  & U \ar[rrdd]^-{b} \ar[rr]^-{(a,b)} & &
  V\times _{\Sym ^{\infty }(X/S)^+}V
  \ar[dd]^-{} \ar[rr]^-{} & & V \ar[dd]^-{G} \\ \\
  & & & V \ar[rr]^-{G} & & \Sym ^{\infty }(X/S)^+
  \enddiagram
  \end{equation}
\fi %%%%%%%%%%%%%%%%%%%%%%%%%%%%%%%%%%%%%%%%%%%%%%%%%%%%%%%%%%%

Consider also the corresponding ``underlying" commutative diagrams
    \begin{equation}
    \label{muhomory2a}
    \xymatrix{
    \Spec (K)^{\hat {\; }}\ar@/_/[dddr]^-{} \ar@/^/[drrr]^-{}
    \ar@{.>}[dr]^-{\hspace{-1mm}\exists ! \hat h} \\
    & \hat W \ar[dd] \ar[rr] & & \hat V \ar[dd]^-{\Delta } \\ \\
    & \hat U \ar[rr]^-{(\hat a,\hat b)} & &
    \hat V\times _{S^{\infty }(X/S)^+}\hat V}
    \end{equation}
and
  \begin{equation}
  \label{muhomory2b}
  \diagram
  \hat V\times _{S^{\infty }(X/S)^+}\hat V
  \ar[dd]^-{} \ar[rr]^-{} & & \hat V \ar[dd]^-{\hat G} \\ \\
  \hat V \ar[rr]^-{\hat G} & & S^{\infty }(X/S)^+
  \enddiagram
  \end{equation}
where $\hat h$ exists and unique due to the commutativities coming from the commutativities in the diagram (\ref{reppa2}).

\iffalse %%%%%%%%%%%%%%%%%%%%%%%%%%%%%%%%%%%%%%%%%%%%%%%%%%%
Summarizing, we obtain the commutative diagram
  \begin{equation}
  \label{muhomory2}
  \diagram
  \Spec (K)^{\hat {\; }} \ar[rrrd]^-{}
  \ar@{.>}[rd]_-{\hspace{+10mm}\exists ! \hat h}
  \ar[rddd]^-{} & & & & & \\
  & \hat W \ar[dd]_-{} \ar[rr]^-{} & & \hat V \ar[rrdd]^-{\id }
  \ar[dd]^-{\Delta } & & \\ \\
  & \hat U \ar[rrdd]^-{\hat b} \ar[rr]^-{(\hat a,\hat b)} & &
  \hat V\times _{S^{\infty }(X/S)^+}\hat V
  \ar[dd]^-{} \ar[rr]^-{} & & \hat V \ar[dd]^-{\hat G} \\ \\
  & & & \hat V \ar[rr]^-{\hat G} & & S^{\infty }(X/S)^+
  \enddiagram
  \end{equation}
\fi %%%%%%%%%%%%%%%%%%%%%%%%%%%%%%%%%%%%%%%%%%%%%%%%%%%%%%%%

Clearly, the commutative diagrams (\ref{muhomory1a}), (\ref{muhomory1b}), (\ref{muhomory2a}) and (\ref{muhomory2b}) can be joined in to one large commutative diagram by means of the morphisms
  $$
  \hat U\to U\; ,\; \; \hat V\to V\; ,\; \; \hbox{etc}
  $$

One of the subdiagrams of that join is the commutative square
  $$
  \diagram
  \hat V\times _{S^{\infty }(X/S)^+}\hat V \ar[dd]_-{} \ar[rr]^-{} & & V \ar[dd]^-{} \\ \\
  \hat V \ar[rr]^-{} & &
  \Sym ^{\infty }(X/S)^+
  \enddiagram
  $$
As we know from the first part of the proof, applied to the case when $U=V$, the diagonal composition
  $$
  \hat V\times _{S^{\infty }(X/S)^+}\hat V
  \to \Sym ^{\infty }(X/S)^+
  $$
is an \'etale neighbourhood of the point $P$.

Since the diagrams
  $$
  \xymatrix{
  \hat U \ar[rr]^-{}
  \ar[ddrr]^-{\hat b} & &
  \hat V\times _{S^{\infty }(X/S)^+}\hat V
  \ar[dd]^-{} \\ \\
  & & \hat V
  }
  $$
and
    $$
    \xymatrix{
    \hat U \ar[dd]^-{} \ar[rr]^-{\hat b} & &
    \hat V \ar[dd]^-{} \\ \\
    U \ar[rr]^-{b} & & V}
    $$
are commutative, we see that the square
  $$
  \xymatrix{
  \hat U \ar[dd]_-{} \ar[rr]^-{} & &
  \hat V\times _{S^{\infty }(X/S)^+}\hat V \ar[dd]^-{} \\ \\
  U \ar[rr]^-{F} & & \Sym ^{\infty }(X/S)^+
  }
  $$
is commutative.

The left vertical arrow in the latter square is \'etale, and the morphism $F$ is \'etale by assumption. Therefore, their composition is \'etale, and we obtain the commutative diagram
  \begin{equation}
  \label{ehma}
  \xymatrix{
  \hat U \ar[rr]^-{}
  \ar[ddrr]^-{} & &
  \hat V\times _{S^{\infty }(X/S)^+}\hat V
  \ar[dd]^-{} \\ \\
  & & \Sym ^{\infty }(X/S)^+
  }
  \end{equation}
in which the morphisms targeted at $\Sym ^{\infty }(X/S)^+$ are \'etale.

Now, if $f:Y\to Y'$ is a morphism between locally Noetherian schemes over a space $\bcZ $, if the structural morphisms $Y\to \bcZ $ and $Y'\to \bcZ $ are \'etale, with regard to the atlas on $\bcZ $, then $f$ is also \'etale. This is an obvious modification of Tag 03FV in \cite{StacksProject}. Applying this property to the diagram (\ref{ehma}), we obtain that the morphism
  $$
  \hat U\to \hat V\times _{S^{\infty }(X/S)^+}\hat V
  $$
is \'etale.

As \'etale morphisms are stable under base change, the Cartesian square from the diagram (\ref{muhomory2a}) then shows that the morphism
  $$
  \hat W\to \hat V
  $$
is \'etale. And since the morphisms $\hat V \to V$ is \'etale, the composition
  $$
  \hat W\to \hat V\to V
  $$
is \'etale. Since $G$ is \'etale by assumption, we see that the composition
  \begin{equation}
  \label{vottak}
  \hat W\to \hat V\to V\stackrel{G}{\lra }\Sym ^{\infty }(X/S)^+
  \end{equation}
is also \'etale.

Finally, analyzing the above join of the commutative diagrams (\ref{muhomory1a}), (\ref{muhomory1b}), (\ref{muhomory2a}) and (\ref{muhomory2b}) by means of the extension morphisms, we see that the composition (\ref{vottak}) is the same as the composition
  $$
  \hat W\to W\to U\stackrel{b}{\lra }
  V\stackrel{G}{\lra }\Sym ^{\infty }(X/S)^+\; .
  $$

Thus, we have obtained the commutative diagram
  $$
  \diagram
  \hat W \ar[dd]_-{} \ar[rr]^-{} & & V \ar[dd]^-{} \\ \\
  U \ar[rr]^-{} & & \Sym ^{\infty }(X/S)^+
  \enddiagram
  $$
\iffalse %%%%%%%%%%%%%%%%%%%%%%%%%%%%%%%%%%%%%%%%%%%%%%%%%%%%
   $$
   \xymatrix{
   & & \hat W \ar[lldd]_-{} \ar[rrdd]^-{} & & \\ \\
   U \ar[rrdd]_-{} & & & & V \ar[lldd]^-{} \\ \\
   & & \Sym ^{\infty }(X/S)^+ & & }
   $$
\fi %%%%%%%%%%%%%%%%%%%%%%%%%%%%%%%%%%%%%%%%%%%%%%%%%%%%%%%%%
whose diagonal composition
  \begin{equation}
  \label{neuzhtokonez}
  \hat W\to \Sym ^{\infty }(X/S)^+
  \end{equation}
is \'etale.

Analyzing the commutative diagrams above, it is easy to see that $P$ factorizes through (\ref{neuzhtokonez}), so that the latter morphism is an \'etale neighbourhood of $P$.
\end{pf}

\subsection{Rational curves on the locally ringed site of $0$-cycles}

Theorem \ref{cofilter} has the following important implication. Namely, since all the items of Tag 00Y5 in \cite{StacksProject} are now satisfied, the stack functor
  $$
  \stalk _P:
  \Shv (\Sym ^{\infty }(X/S)^+_{\Nis \mhyphen \et })\to
  \Sets
  $$
induces a point of the topos 
  $$
  \Shv (\Sym ^{\infty }(X/S)^+_{\Nis \mhyphen \et })
  $$ 
by Lemma 7.31.7 in loc.cit. In particular, we obtain the full-fledged stalk
  $$
  \bcO _{\Sym ^{\infty }(X/S)^+\! ,\, P}=
  \stalk _P\, (\bcO _{\Sym ^{\infty }(X/S)^+})
  $$

Moreover, the ringed site 
  $$
  \Sym ^{\infty }(X/S)^+_{\Nis \mhyphen \et }
  $$ 
is a locally ringed site in the sense of the definition appearing in Exercise 13.9 on page 512 in \cite{SGA4-1} (see page 313 in the newly typeset version), as well as in the sense of a sightly different Tag 04EU in \cite{StacksProject}. This is explained in Section \ref{kaehler}.

For short of notation, let us write
  $$
  \bcO _P=\bcO _{\Sym ^{\infty }(X/S)^+\! ,\, P}\; .
  $$
This should not lead to a confusion, as the point $P$ is a point on $\Sym ^{\infty }(X/S)^+$. Since the site $\Sym ^{\infty }(X/S)^+_{\Nis \mhyphen \et }$ is locally ringed, for each point $P$ on this site the stalk $\bcO _P$ is a local ring by the same Tag 04ET in \cite{StacksProject}. Then we also have the maximal ideal
  $$
  \gom _P\subset \bcO _P
  $$
and the residue field
  $$
  \kappa (P)=\bcO _P/\gom _P
  $$
at the point $P$.

The stalk functor also gives us the stalks
  $$
  \Omega ^1_{\Sym ^{\infty }(X/S)^+\!,\, P}=
  \stalk _P\, (\Omega ^1_{\Sym ^{\infty }(X/S)^+})
  $$
and
  $$
  T_{\Sym ^{\infty }(X/S)^+\! ,\, P}=
  \stalk _P\, (T_{\Sym ^{\infty }(X/S)^+})
  $$
at $P$. Tensoring by $\kappa (P)$ we obtain the vector spaces
  $$
  \Omega ^1(P)=\Omega ^1_{\Sym ^{\infty }(X/S)^+}(P)=
  \Omega ^1_{\Sym ^{\infty }(X/S)^+\!,\, P}
  \otimes _{\bcO _P}\kappa (P)
  $$
and
  $$
  T(P)=T_{\Sym ^{\infty }(X/S)^+}(P)=
  T_{\Sym ^{\infty }(X/S)^+\! ,\, P}
  \otimes _{\bcO _P}\kappa (P)
  $$
over the residue field $\kappa (P)$.

The second vector space $T(P)$ is then our {\it tangent space} to the space of $0$-cycles $\Sym ^{\infty }(X/S)^+$ at the point $P$. Notice that, since $\Sym ^{\infty }(X/S)^+$ is an abelian group object in the category of Nisnevich sheaves on locally Noetherian schemes over $S$, whenever $S$ is the spectrum of a field $k$, all tangent spaces $T(P)$ at $k$-rational points $P$ are uniquely determined by the tangent space $T(0)$ at the zero point $0$ on $\Sym ^{\infty }(X/S)^+$ provided by the section of the structural morphism from $X$ to $S$. In other words, one can develop a Lie theory on $\Sym ^{\infty }(X/S)^+$.

Now we are fully equipped to promote the idea of understanding of rational equivalence of $0$-cycles as rational connectivity on the space $\Sym ^{\infty }(X/S)^+$. First of all, looking at any scheme $U$ over $S$ as a representable sheaf, we have the corresponding locally ringed site $U_{\Nis \mhyphen \et }$. Then a {\it regular morphism} from $U$ to $\Sym ^{\infty }(X/S)^+_{\Nis \mhyphen \et }$ is just a morphism of locally ringed sites
  $$
  U_{\Nis \mhyphen \et }\to
  \Sym ^{\infty }(X/S)^+_{\Nis \mhyphen \et }
  $$
in the sense of Tag 04HA in \cite{StacksProject}. Notice that since \'etale morphisms are stable under base change, if $U\to \Sym ^{\infty }(X/S)^+$ is a morphism of sheaves, then it induces the corresponding morphism of locally ringed sites.\label{proverit' !!!}

A rational curve on $\Sym ^{\infty }(X/S)^+$ is a morphism of sheaves
  $$
  f:\PR ^1\to \Sym ^{\infty }(X/S)^+\; .
  $$
If
  $$
  P:\Spec (K)\to \Sym ^{\infty }(X/S)^+
  $$
is a point on the sheaf $\Sym ^{\infty }(X/S)^+$, then we will be saying that $f$ passes through the point $P$ if $P$, as a morphism to $\Sym ^{\infty }(X/S)^+$, factorizes through the morphism $f:\PR ^1\to \Sym ^{\infty }(X/S)^+$.

Now, two points $P$ and $Q$ on $\Sym ^{\infty }(X/S)^+$ are {\it elementary rationally connected} if there exists a rational curve on $\Sym ^{\infty }(X/S)^+$ passing through $P$ and $Q$. The points $P$ and $Q$ are said to be {\it rationally connected} if there exists a finite set of points $R_1,\ldots ,R_n$ on $\Sym ^{\infty }(X/S)^+$, such that $R_1=P$, $R_n=Q$ and $R_i$ is elementary rationally connected to $R_{i+1}$ for each $i\in \{ 1,\ldots ,n-1\} $. If any two points on $\Sym ^{\infty }(X/S)^+$ are rationally connected, then we will say that this space is rationally connected.

Let
  $$
  P:\Spec (K)\to \Sym ^{\infty }(X/S)^+
  \qqand
  Q:\Spec (L)\to \Sym ^{\infty }(X/S)^+
  $$
be two points on $\Sym ^{\infty }(X/S)^+$, represented by morphisms from the spectra of two fields $K$ and $L$ respectively. Suppose, in addition, that the fields $K$ and $L$ are embedded in to a common field, in which case we can replace both $K$ and $L$ by their composite $KL$. Then we can assume, without loss of generality, that $K=L$. In such a case, the points $P$ and $Q$, as morphisms from the scheme $\Spec (K)$ to the sheaf $\Sym ^{\infty }(X/S)^+$ induce two sections $s_P$ and $s_Q$ in
  $$
  \Sym ^{\infty }(X/S)^+(\Spec (K))=
  \bcZ _0^{\infty }(X/S)(\Spec (K))=
  $$
  $$
  z(X/S,0)_{\infty }(\Spec (K))\; .
  $$
Assume, in addition, that
  $$
  S=\Spec (K)\; .
  $$
Then $s_P$ and $s_Q$, as elements of the group
  $$
  z(X/\Spec (K),0)_{\infty }(\Spec (K))\; ,
  $$
are two $0$-cycles on the scheme $X$ over $\Spec (K)$. And since relative $0$-cycles are representable, see Section \ref{relcycles}, rational connectivity of the points $P$ and $Q$ on $\Sym ^{\infty }(X/S)^+$ is equivalent to rational equivalence of the $0$-cycles $s_P$ and $s_Q$ on the scheme $X$.
This all means that we can look at rational connectedness between points on $\Sym ^{\infty }(X/S)^+$ as the generalized rational equivalence in the relative setting.

Let $V$ be an arbitrary smooth projective variety over $k$. According to Koll\'ar, \cite{KollarRatCurvesOnVar}, if we wish to show that $V$ is rationally connected, we should to two steps. The first one is that we need to find a rational curve
  $$
  f:\PR ^1\to V
  $$
on $V$. If the first step is done, then we need to show that the rational curve $f$ is free on $V$, i.e that the numbers
  $$
  a_1\geq \ldots \geq a_n
  $$
in the decomposition
  $$
  f^*T_V=
  \bcO _{\PR ^1}(a_1)\oplus \ldots \oplus \bcO _{\PR ^1}(a_n)
  $$
have appropriate positivity, where $T_V$ is the tangent sheaf on the variety $V$, see Section II.3 in the canonical book \cite{KollarRatCurvesOnVar}, or many other sources about free curves on varieties.

Now, since we have the tangent sheaf $T_{\Sym ^{\infty }(X/k)^+}$ for our surface $X$ over $k$, we can try to do the same on the space $\Sym ^{\infty }(X/k)^+$. Namely, we should first find a rational curve
  $$
  f:\PR ^1\to \Sym ^{\infty }(X/k)^+
  $$
on the space of $0$-cycles. 

Of course, we do not know (at the moment) whether the tangent sheaf 
  $$
  T_{\Sym ^{\infty }(X/k)^+}
  $$ 
is locally free on the site $\Sym ^{\infty }(X/S)^+_{\Nis \mhyphen \et }$, and, accordingly, we do not know whether the pullback 
  $$
  f^*T_{\Sym ^{\infty }(X/k)^+}
  $$ 
decomposes in to the direct sum of Serre twists. 

But it is not hard to show that this pullback is a coherent sheaf on the projective line $\PR ^1$ over $k$. Being a coherent sheaf, it decomposes uniquely in to a direct sum of a torsion sheaf and a locally free sheaf, see, for example, Proposition 5.4.2. in \cite{ChenKrause}. 

Then
  $$
  f^*T_{\Sym ^{\infty }(X/k)^+}=
  \bcO _{\PR ^1}(a_1)\oplus \ldots \oplus \bcO _{\PR ^1}(a_n)
  \oplus \bcT \; ,
  $$
where $\bcT $ is a sheaf on $\PR ^1$ to be understood. 

Though the sheaf $\bcT $ is possibly non-zero and mysterious, we still can try to apply the same line of arguments as in the proof of Theorem 3.7 in \cite{KollarRatCurvesOnVar} or Proposition 4.8 in \cite{Debarre}.

\bigskip

\begin{appendix}

\section{Categorical monoids and group completions}
\label{freemongrcompl}

\subsection{Free monoids in cartesian categories}

Let $\catS $ be a cartesian monoidal category, so that the terminal object $\! \term \! $ is the monoidal unit in $\catS $. Denote by $\Mon (\catS )$ the full subcategory of monoids\footnote{all monoids in this paper will be commutative by default}, and by $\Ab (\catS )$ the full subcategory of abelian group objects in the category $\catS $. Assume that $\catS $ is closed under finite colimits and countable coproducts which are distributive with regard to the cartesian product in $\catS $. Then the forgetful functor from $\Mon (\catS )$ to $\catS $ has left adjoint which can be constructed as follows.

For any object $\bcX $ in $\catS $ and for any natural number $d$ let $\bcX ^d$ be the $d$-fold monoidal product of $\bcX $. Consider the $d$-th symmetric power
  $$
  \Sym ^d(\bcX )\; ,
  $$
i.e. the quotient of the object $\bcX ^d$ by the natural action of the $d$-th symmetric group $\Sigma _d$ in the category $\catS $. In particular,
  $$
  \Sym ^0(\bcX )=\term
  \qqand
  \Sym ^1(\bcX )=\bcX \; .
  $$
The coproduct
  $$
  \coprod _{d=0}^{\infty }\Sym ^d(\bcX )
  $$
is a monoid, whose concatenation product
  $$
  \coprod _{d=0}^{\infty }\Sym ^d(\bcX )\times
  \coprod _{d=0}^{\infty }\Sym ^d(\bcX )\to
  \coprod _{d=0}^{\infty }\Sym ^d(\bcX )
  $$
is induced by the obvious morphism
  $$
  \coprod _{d=0}^{\infty }\bcX ^{d}\times
  \coprod _{d=0}^{\infty }\bcX ^{d}\to
  \coprod _{d=0}^{\infty }\bcX ^{d}
  $$
and the embeddings of $\Sigma _i\times \Sigma _j$ in to $\Sigma _{i+j}$. The unit
  $$
  \term \to \coprod _{d=0}^{\infty }\Sym ^d(\bcX )
  $$
identifies $\term $ with $\bcX ^{(0)}$. This monoid will be called the {\it free monoid} generated by $\bcX $ and denoted by $\NN (\bcX )$. Thus,
  $$
  \NN (\bcX )=\coprod _{d=0}^{\infty }\Sym ^d(\bcX )\; .
  $$
For example,
  $$
  \NN (\term )=\NN \; .
  $$
It is easy to verify that the functor
  $$
  \NN :\catS \to \Mon (\catS )
  $$
is left adjoint to the forgetful functor from $\Mon (\catS )$ to $\catS $.

The full embedding of $\Ab (\catS )$ in to $\Mon (\catS )$ admits left adjoint, if we impose some extra assumption on the category $\catS $. Namely, let $\bcX $ be a monoid in $\catS $, and look at the obvious diagonal morphism
  \begin{equation}
  \label{diag}
  \Delta :\bcX \to \bcX \times \bcX
  \end{equation}
in the category $\catS $, which is also a morphism in the category $\Mon (\catS )$. The terminal object $*$ in the category $\catS $ is a trivial monoid, i.e. a terminal object in the category $\Mon (\catS )$.

Assume there exists a co-Cartesian square
  \begin{equation}
  \label{completiondiagr}
  \diagram
  \bcX \ar[rr]^-{\Delta }
  \ar[dd]^-{}
  & & \bcX \times \bcX \ar[dd]^-{} \\ \\
  \term \ar[rr]^-{} & & \bcX ^+
  \enddiagram
  \end{equation}
in the category of monoids $\Mon (\catS )$. Then $\bcX ^+$ is an abelian group object in the category $\catS $.

Let
  $$
  \iota _{\bcX }:\bcX \to \bcX ^+
  $$
be the composition of the canonical embedding
  $$
  \iota _1:\bcX \to \bcX \times \bcX \; ,
  $$
  $$
  x\mapsto (x,0)
  $$
with the projection
  $$
  \pi _{\bcX }:\bcX \times \bcX \to \bcX ^+\; .
  $$
If
  $$
  f:\bcX \to \bcY
  $$
is a morphism of monoids and $\bcY $ is an abelian group object in $\catS $, the precomposition of the homomorphism
  $$
  (f,-f):\bcX \times \bcX \to \bcY \; ,
  $$
sending $(x_1,x_2)$ to $f(x_1)-f(x_2)$ with the diagonal embedding is $0$, whence there exists a unique group homomorphism $h$ making the diagram
  $$
  \xymatrix{
  \bcX \ar[rr]^-{\iota _{\bcX }}
  \ar[ddrr]_-{f} & &
  \bcX ^+ \ar@{.>}[dd]^-{\hspace{+1mm}\exists ! h} \\ \\
  & & \bcY
  }
  $$
commutative.

This all shows that $\bcX ^+$ is nothing else but the the {\it group completion} of the monoid $\bcX $, and the group completion functor
  $$
  -^+:\Mon (\catS )\to \Ab (\catS )
  $$
is left adjoint to the forgetful functor from $\Ab (\catS )$ to $\Mon (\catS )$.

For example,
  $$
  \ZZ =\NN ^+
  $$
is the group completion of the free monoid $\NN $, generated by the terminal object $\term $ in the category $\catS $.

Notice that, as the categories $\Mon (\catS )$ and $\Ab (\catS )$ are pointed, one can show the existence of the canonical isomorphism of monoids
  $$
  (\bcX \times \bcX )^+\stackrel{\sim }{\to }
  \bcX ^+\times \bcX ^+\; .
  $$
In other words, the group completion functor is monoidal.

It is useful to understand how all these constructions work for set-theoretical monoids. Since monoids are not groups, some care is in place here.

Let $M$ be a monoid in the category of sets $\Sets $, written additively, and assume first that we are given with a submonoid $N$ in $M$. To understand what would be the quotient monoid of $M$ by $N$, we define a relation
  $$
  R\subset M\times M
  $$
saying that, for any two elements $m,m'\in M$,
  \begin{equation}
  \label{defofeq}
  mRm'\; \Leftrightarrow \; \exists n,n'\in N\; \hbox{with}\;
  m+n=m'+n'\; .
  \end{equation}
Then $R$ is a congruence relation on $M$, i.e. an equivalence relation compatible with the operation in $M$. Indeed, the reflexivity and symmetry are obvious. Suppose that we have three elements
  $$
  m,m',m''\in M\; ,
  $$
and
  $$
  \exists n,n'\in N,\; \hbox{such that}\; m+n=m'+n'\; .
  $$
and
  $$
  \exists l',l''\in N,\; \hbox{such that}\;
  m'+l'=m''+l''\; .
  $$
Then
  $$
  m+n+l'=m'+n'+l'=m''+l''+n'\; .
  $$
Clearly,
  $$
  n+l',\; l''+n'\in N\; ,
  $$
and we get transitivity. Thus, $R$ is an equivalence relation.

Let $M/N$ be the corresponding quotient set, and let
  $$
  \pi :M\to M/N
  $$
  $$
  m\mapsto [m]
  $$
be the quotient map. The structure of a monoid on $M/N$ is obvious,
  $$
  [m]+[\tilde m]=[m+\tilde m]\; ,
  $$
and since $M$ is a commutative monoid\footnote{recall that, within this paper, all monids are commutative by defaul}, it follows easily that the map $\pi $ is a homomorphism of monoids. In other terms, the above relation $R$ on $M$ is a congruence relation.

Moreover, the quotient homomorphism
  $$
  M\to M/N
  $$
enjoys the standard universal property, loc.cit. To be more precise, for any homomorphism of monoids
  $$
  f:M\to T\; ,
  $$
such that
  $$
  N\subset \ker (f)=\{ m\in M\; |\; f(m)=0\} \; ,
  $$
there exists a commutative diagram of type
  $$
  \xymatrix{
  M \ar[rr]^-{} \ar[ddrr]_-{} & &
  M/N \ar@{.>}[dd]^-{\hspace{+1mm}\exists !} \\ \\
  & & T
  }
  $$

Now, let $M\times M$ be the product monoid, let
  $$
  \Delta :M\to M\times M
  $$
be the diagonal homomorphism, and let
  $$
  \Delta (M)
  $$
be the set-theoretical image of the homomorphism $\Delta $. Trivially, $\Delta (M)$ is a submonoid in the product monoid $M\times M$, and we can construct the quotient monoid
  $$
  M^+=(M\times M)/\Delta (M)\; ,
  $$
using the procedure explained above. The universal property of the quotient monoid gives us that the diagram
  \begin{equation}
  \label{mohoviki}
  \diagram
  M\ar[rr]^-{\Delta } \ar[dd]^-{} & &
  M\times M \ar[dd]^-{} \\ \\
  \term \ar[rr]^-{} & & M^+
  \enddiagram
  \end{equation}
is pushout in the category $\Mon (\Sets )$. It follows that $M^+$ is the group completion of $M$ in the sense of our definition given for the general category $\catS $.

Clearly, the composition
  $$
  \xymatrix{
  M \ar[rr]^-{m\mapsto (m,0)} \ar[ddrr]_-{\iota _M} & &
  M\times M \ar[dd]^-{} \\ \\
  & & M^+
  }
  $$
is a homomorphism of monoids. If
  $$
  f:M\to A
  $$
is a homomorphism from the monoid $M$ to an abelian group $A$, then we define a homomorphism of monoids
  $$
  M\times M\to A
  $$
sending
  $$
  (m_1,m_2)\mapsto f(m_1)-f(m_2)\; ,
  $$
and the universal property of the diagram (\ref{mohoviki}) gives us the needed commutative diagram
  $$
  \xymatrix{
  M \ar[rr]^-{\iota _M} \ar[ddrr]_-{} & &
  M^+ \ar@{.>}[dd]^-{\hspace{+1mm}\exists !} \\ \\
  & & A
  }
  $$

Moreover, if $M$ is cancellative, the diagram (\ref{mohoviki}) is not only a pushout square in $\Mon (\Sets )$ but also a pullback square in $\Sets $.

Indeed, if
  $$
  (m_1,m_2),\; (m_1',m_2')\in M\times M\; ,
  $$
then, according to (\ref{defofeq}),
  $$
  \exists n,n'\in M
  $$
such that
  $$
  (m_1,m_2)+(n,n)=(m_1',m_2')+(n',n')
  $$
in $M\times M$, or, equivalently,
  \begin{equation}
  \label{kott}
  m_1+n=m_1'+n'
  \qqand
  m_2+n=m_2'+n'\; .
  \end{equation}

Now, suppose we want to find $h$ completing a commutative diagram of type
    \begin{equation}
    \label{ryzhiikot}
    \xymatrix{
    T\ar@/_/[dddr] \ar@/^/[drrr]^-{f}
    \ar@{.>}[dr]^-{\hspace{-1mm}\exists ! h} \\
    & M \ar[dd] \ar[rr] & & M\times M \ar[dd]^-{\pi } \\ \\
    & \term \ar[rr] & & M^+}
    \end{equation}
in the category $\Sets $. If $(m_1,m_2)$ is an element of $M\times M$, the equivalence class $[m_1,m_2]$ is $0$ in $M^+$, i.e. the ordered pair $(m_1,m_2)$ is equivalent to $(0,0)$ in $M\times M$ modulo the subtractive submonoid $\Delta (M)$, if and only if, by (\ref{kott}),
  $$
  m_1+n=n'
  \qqand
  m_2+n=n'\; ,
  $$
whence
  $$
  m_1+n=m_2+n\; .
  $$
Since $M$ is a cancellation monoid, the latter equality gives us that $m_1=m_2$, i.e. $(m_1,m_2)$ is in $\Delta (M)$. In other words, $[m_1,m_2]=0$ in $M^+$ if and only if $(m_1,m_2)$ is in $\Delta (M)$. And as the diagram (\ref{ryzhiikot}) is commutative without $h$, it follows that the set-theoretical image of the map $f$ is in $\Delta (M)$. It follows that $f$ factorizes through $\Delta $, i.e. the needed map $h$ exists.

Thus, we see that the abstract constructions relevant to group completions are generalizations of the standard constructions in terms of set-theoretical monoids.

All the same arguments apply when $\catS $ is the category $\PShv (\catC )$ of set valued presheaves on a category $\catC $, as all limits and colimits in $\PShv (\catC )$ are sectionwise. Thus, for any monoid $\bcX $ in $\PShv (\catC )$ the group completion $\bcX ^+$ exists and it is a section wise group completion. If $\bcX $ is cancellative, and this is equivalent to saying that $\bcX $ is section wise cancellative, then the diagram (\ref{completiondiagr}) is Cartesian in $\PShv (\catC )$.

\subsection{Connective v.s. disconnective monoids}
\label{condiscon}

Now let us come back to the general setting. Let again $\bcX $ be a monoid in $\catS $. The notion of a cancellative monoid can be categorified as follows. A morphism
  $$
  \iota :\NN \to \bcX
  $$
in the category $\Mon (\catS )$, that is a homomorphism of monoids from $\NN $ to $\bcX $, is uniquely defined by the restriction
  $$
  \iota (1):\term \to \bcX
  $$
of $\alpha $ on to the subobject $\term =\Sym ^1(\term )$ of the object 
  $$
  \NN =\coprod _{d=0}^{\infty }\Sym ^d(\term )
  $$ 
in $\catS $. 

Vice versa, as soon as we have a morphism $\term \to \bcX $ in the category $\catS $, it uniquely defines the obvious morphism $\iota :\NN \to \bcX $ in the category $\Mon (\catS )$. The homomorphism of monoids $\iota $ will be said to be {\it cancellative} if the composition
  $$
  \add _{\iota (1)}:\bcX \simeq \bcX \times *\stackrel{\id \times \iota (1)}{\lra }\bcX \times \bcX \to \bcX \; ,
  $$
that is the addition of $\iota (1)$ on $\bcX $, is a monomorphism in $\catS $. 

The monoid $\bcX $ is a {\it cancellation monoiod} if any homomorphism $\iota :\NN \to \bcX $ is cancellative.

Clearly, if $\bcX $ is a monoid in $\PShv (\catC )$, then $\bcX $ is cancellative if and only if it is section wise cancellative.

A {\it pointed monoid} in $\catS $ is a pair $(\bcX ,\iota )$, where $\bcX $ is a monoid in $\catS $ and $\iota $ is a morphism of monoids from $\NN $ to $\bcX $. 

A {\it graded pointed monoid} is a triple $(\bcX ,\iota ,\sigma )$, where $(\bcX ,\iota )$ is a pointed monoid and $\sigma $ is a morphism of monoids from $\bcX $ to $\NN $, such that
  $$
  \sigma \circ \iota =\id _{\NN }\; .
  $$

If $\bcX $ is a pointed graded monoid in $\catS $, for any natural number $d\in \NN $ one can consider the cartesian square
  $$
  \diagram
  \bcX _d\ar[rr]^-{} \ar[dd]^-{}
  & & * \ar[dd]^-{d} \\ \\
  \bcX \ar[rr]^-{\sigma } & & \NN
  \enddiagram
  $$
in the category $\catS $. The addition of $\iota (1)$ in $\bcX $ induces morphisms
  $$
  \bcX _d\to \bcX _{d+1}
  $$
for all $d\geq 0$. Let $\bcX _{\infty }$ be the colimit
  $$
  \bcX _{\infty }=
  \colim (\bcX _0\to \bcX _1\to \bcX _2\to \dots )
  $$
in $\catS $. Equivalently, $\bcX _{\infty }$ is the coequalizer of the addition of $\iota (1)$ in $\bcX $ and the identity automorphism of $\bcX $. Since filtered colimits commute with finite products, there is a canonical isomorphism between the colimit of the obvious diagram composed by the objects $\bcX _d\times \bcX _{d'}$, for all $d,d'\geq 0$, and the product $\bcX _{\infty }\times \bcX _{\infty }$. Since the colimit of that diagram is the colimit of its diagonal, this gives the canonical morphism from $\bcX _{\infty }\times \bcX _{\infty }$ to $\bcX _{\infty }$. The latter defines the structure of a monoid on $\bcX _{\infty }$, such that the canonical morphism
  $$
  \pi :\bcX =\coprod _{d\geq 0}\bcX _d\to \bcX _{\infty }
  $$
is a homomorphism of monoids in $\catS $. We call $\bcX _{\infty }$ the {\it connective} monoid associated to the pointed graded monoid $\bcX $.

Notice that if the category $\catS $ is exhaustive\footnote{see {{\tt https://ncatlab.org/nlab/show/exhaustive+category}}}, monomorphicity of the morphisms $\bcX _d\to \bcX _{d+1}$ yields that the transfinite compositions $\bcX _d\to \bcX _{\infty }$ are monomorphisms too. The morphisms $\bcX _d\to \bcX _{d+1}$ are monomorphic, for example, if $\bcX $ is a cancelation monoid.

Now assume that the colimit $\bcX _{\infty }^+$ exists in the category $\Mon (\catS )$. Since $\bcX _{\infty }$ is the coequalizer of $\add _{\iota (1)}$ and $\id _{\bcX }$, the group completion $\bcX _{\infty }^+$ is the coequalizer of the corresponding homomorphism $\add _{\iota (1)}^+:\bcX ^+\to \bcX ^+$ and $\id _{\bcX ^+}$. It follows that the sequence
  $$
  0\to \ZZ \stackrel{\iota ^+}{\lra }\bcX ^+\to
  \bcX _{\infty }^+\to 0
  $$
is short exact. Moreover, this sequens splits by the morphism $\sigma ^+$. This gives us that
  $$
  \bcX ^+=\ZZ \oplus \bcX _{\infty }^+
  $$
in the abelian category $\Ab (\catS )$.

A typical example of a pointed graded monoid in $\catS $ is the free monoid
  $$
  \NN (\bcX )=\coprod _{d=0}^{\infty }\Sym ^d(\bcX )\; ,
  $$
where $\bcX $ is a pointed object in $\catS $, i.e. the morphism from $\bcX $ to the terminal object $*$ has a section. For this pointed graded monoid we have that
  $$
  \NN (\bcX )_d=\Sym ^d(\bcX )\; ,
  $$
for all natural numbers $d$, and the pointing of each symmetric power $\Sym ^d(\bcX )$ is induced by the pointing of $\bcX $ in the obvious way. The section gives embeddings
  $$
  \Sym ^d(\bcX )\to \Sym ^{d+1}(\bcX )\; ,
  $$
and the corresponding connective monoid
  $$
  \NN (\bcX )_{\infty }=\colim _d\; \Sym ^d(\bcX )
  $$
will be denoted by $\Sym ^{\infty }(\bcX )$ and called the {\it free connective monoid} of the object $\bcX $. Then, of course,
  $$
  \Sym ^{\infty }(\bcX )^+=\NN (\bcX )_{\infty }^+\; .
  $$
Moreover, both free monoids $\NN (\bcX )$ and $\NN (\bcX )_{\infty }$ are cancellative monoids in $\catS $.

\subsection{In sheaves and presheaves}
\label{sheavespresheaves}

Now, let $\catC $ be a cartesian monoidal category with a terminal object $\! \term \, $, closed under finite fibred products and equipped with a subcanonical topology $\tau $. Let $\PShv (\catC )$ be the category of set valued presheaves on $\catC $, and let $\Shv (\catC _{\tau })$ be the full subcategory in $\PShv (\catC )$ of sheaves on $\catC $ with regard to the topology $\tau $. Since the category $\catC $ is cartesian, so are the categories $\PShv (\catC )$ and $\Shv (\catC _{\tau })$, and therefore we can consider the monoids in the categories of sheaves and pre-sheaves. Our aim is now to apply the constructions above in the case when
  $$
  \catS =\PShv (\catC )
  \qor
  \catS =\Shv (\catC _{\tau })\; .
  $$

The Yoneda embedding
  $$
  h:\catC \to \PShv (\catC )
  $$
is a continuous functor, i.e. it preserves limits. It follows that, if $\bcX $ is a monoid in $\PShv (\catC )$, then it is equivalent to saying that $\bcX $ is a section wise monoid. Moreover, the diagonal morphism (\ref{diag}) for a presheaf $\bcX $ is diagonal section-wise. It follows that the colimit diagram (\ref{completiondiagr}) exists in $\Mon (\PShv (\catC ))$ and, accordingly, the group completion $\bcX ^+$ is then the section wise group completion of $\bcX $. In particular, the group completion $\bcX ^+$ of the presheaf monoid $\bcX $ is topology free.

Since the sheafification functor $-^{\shf }$ is left adjoint to the forgetful functor from sheaves to presheaves, the latter is right adjoint, and hence it commutes with limits. In particular, the forgetful functor from sheaves to presheaves commutes with products. 

Next, it is well-known that the functor $-^{\shf }$ is exact too, and hence it commutes with products. It follows that $-^{\shf }$ takes monoids to monoids, and abelian groups to abelian groups.

The functor $\NN $ exists for set valued presheave monoids and it is given section wise. Moreover, as we mentioned above, the group completion functor exists for set valued presheaf monoids, and it is also given section wise. It follows that the functors $\NN $ and $-^+$ exist also for sheaves on the site $\catC _{\tau }$, and can be constructed by means of composing of the corresponding functors for presheaves with the sheafification functor.

To be more precise, since the sheafification $-^{\shf }$ is left adjoint, it also commutes with all colimits. And as the functors $\NN $ and $-^+$ are constructed merely by means of products and colimits, we conclude that that these two functors are preserved by sheafification. 

It means that, if $\bcX $ is a set valued sheaf on $\catC _{\tau }$, then, in order to construct the free monoid $\NN (\bcX )$ in the category of sheaves, we first forget the sheaf property on $\bcX $ and construct $\NN (\bcX )$ in the category of pre-sheaves, and then sheafify to get an object in monoids of sheaves on $\catC _{\tau }$. If $\bcX $ is a set valued sheaf monoid, then, in order to construct its group completion in the category of sheaves, we forget the sheaf property on $\bcX $ and construct $\bcX ^+$ in the category of presheaves, and then sheafify to get an object in sheaves. 

Similarly, if $\bcX $ is a pointed graded monoid in presheaves, then it is a pointed graded monoid section wise. The construction of the connective monoid $\bcX _{\infty }$, as an object in the category $\Mon (\PShv (\catC ))$, is then section wise and topology free. But if $\bcX $ is a pointed graded monoid in sheaves, the construction of $\bcX _{\infty }$ follows the rule above. Namely, we first forget the sheaf property of $\bcX $ and construct $\bcX _{\infty }$ section wise, and then sheafify.

As in the previous section, for simplicity of notation, we will write $X$ instead of the sheaf $h_X$, for any object $X$ in $\catC $, and denote objects in $\PShv (\catC )$ and $\Shv (\catC _{\tau })$ by calligraphic letters $\bcX $, $\bcY $, etc.

Notice also that if $X$ is a pointed object of $\catC $ and for any $d$ the $d$-th symmetric power $\Sym ^d(X)$ exists already in $\catC $, then $\NN (X)_{\infty }$ is an $\ind $-object of $\catC $. Recall that an $\ind $-object in $\catC $ is the colimit of the composition of a functor
  $$
  I\to \catC
  $$
with the embedding of $\catC $ in to $\PShv (\catC )$, taken in the category $\PShv (\catC )$, such that the category $I$ is filtered. Such a colimit is section-wise. Since $\catC $ is equipped with a topology, one can also give the definition of a sheaf-theoretical $\ind $-object. An $\ind $-object in $\catC _{\tau }$ is the colimit of the same composition, but now taken in the category $\Shv (\catC _{\tau })$. The latter is obviously the sheafification of the previous $\ind $-object, and therefore it depends on the topology $\tau $. Let $\Ind (\catC )$ be the full subcategory in $\PShv (\catC )$ of $\ind $-objects in $\catC $, and let $\Ind (\catC _{\tau })$ be the full subcategory in $\Shv (\catC _{\tau })$ of $\ind $-objects of $\catC _{\tau }$.

Let us apply these abstract constructions to the case when
  $$
  \catC =\Noe /S
  $$
and
  $$
  \tau =\Nis \; .
  $$
The choice of the topology will be explained in the next section. Now we need to recall relative symmetric powers of locally Noetherian schemes $X$ over $S$.

Assume that the structural morphism
  $$
  X\to S
  $$
satisfies the following property:

\begin{itemize}

\item[(AF)]{}
for any point $s\in S$ and for any finite collection $\{ x_1,\ldots ,x_l\} $ of points in the fibre $X_s$ of the structural morphism $X\to S$ at $s$ there exists a Zariski open subset $U$ in $X$, such that
  $$
  \{ x_1,\ldots ,x_l\} \subset U
  $$
and the composition
  $$
  U\to X\to S
  $$
is a quasi-affine morphism of schemes.

\end{itemize}

Quasi-affine morphisms possess various nice properties, see \cite{StacksProject}, which can be used to prove that if $U\to S$ is a morphism of locally Noetherian schemes and $X$ is AF over $S$ then $X\times _SU$ is AF over $U$. If, moreover, $U$ is AF over $S$ the $X\times _SU$ is AF over $S$.

The property AF is satisfied if, for example, $X\to S$ is a quasi-affine or quasi-projective morphism of schemes, see Prop. (A.1.3) in Paper I in \cite{RydhThesis}.

As we now assume that AF holds true for $X$ over $S$, the $d$-th symmetric group $\Sigma _d$ acts admissibly on the $d$-th fibred product
  $$
  (X/S)^d=X\times _S\ldots \times _SX
  $$
over $S$ in the sense of \cite{SGA1}, Expos\'e V, and the relative symmetric power
  $$
  \Sym ^d(X/S)
  $$
exists in the category $\Noe /S$.

Then, according to the abstract constructions above, we obtain the free monoid $\NN (X/S)$ generated by the scheme $X$ over $S$ in the category $\Shv ((\Noe /S)_{\Nis })$. For every integer $d\geq 0$ the object $\NN (X/S)_d$ is the relative $d$-th symmetric power $\Sym ^d(X/S)$ of $X$ over $S$, and as such it is an object of the category $\Noe /S$. The free monoid of the scheme $X$ over $S$ is nothing else but the coproduct
  $$
  \NN (X/S)=\coprod _{d=0}^{\infty }\Sym ^d(X/S)
  $$
taken in the category $\Shv ((\Noe /S)_{\Nis })$.

Assume, in addition, that the structural morphism $X\to S$ has a section
  $$
  S\to X\; .
  $$
Notice that the terminal object $*$ in the category $\Noe /S$ is the identity morphism of the scheme $S$, and therefore the splitting of the structural morphism $X\to S$ by the section $S\to X$ induces the splitting of $\NN $ from $\NN (X/S)$ in the category $\Shv ((\Noe /S)_{\Nis })$. 

The corresponding connective monoid
  $$
  \Sym ^{\infty }(X/S)=\NN (X/S)_{\infty }=
  \colim _d\, \Sym ^d(X/S)
  $$
is an $\ind $-scheme over $S$. As such it can be considered as an object of the category $\Ind ((\Noe /S)_{\Nis })$.

The colimit
  $$
  \Sym ^{\infty }(X/S)^+
  $$
in the category of monoids in $\Shv ((\Noe /S)_{\Nis })$ is the group completion of the monoid $\Sym ^{\infty }(X/S)$, and, according to what we discussed above, this colimit is nothing but the Nisnevich sheafification of the corresponding section wise colimit.

\end{appendix}

\bigskip

\begin{small}

\end{small}

\bigskip

\begin{small}

\end{small}

\bigskip

\begin{small}

{\sc Department of Mathematical Sciences, University of Liverpool, Peach Street, Liverpool L69 7ZL, England, UK}

\end{small}

\bigskip

\begin{footnotesize}

{\it E-mail address}: {\tt vladimir.guletskii@liverpool.ac.uk}

\end{footnotesize}

\bigskip

\end{document}